\documentclass[11pt, a4paper, oneside, reqno]{amsart}

\usepackage[usenames, dvipsnames]{color}
\definecolor{darkblue}{rgb}{0.0, 0.0, 0.45}
\definecolor{lightblue}{RGB}{240,248,255}
\definecolor{lightblue2}{rgb}{0.68, 0.85, 0.9}
\definecolor{lightcyan}{rgb}{0.88, 1.0, 1.0}
\definecolor{palepink}{rgb}{0.98, 0.85, 0.87}

\usepackage[colorlinks	= true,
raiselinks	= true,
linkcolor	= darkblue, 
citecolor	= Mahogany,
urlcolor	= ForestGreen,
pdfauthor	= {Peyman Mohajerin Esfahani},
pdftitle	= {},
pdfkeywords	= {},
pdfsubject	= {},
plainpages	= false]{hyperref}
\usepackage{dsfont,amssymb,amsmath,enumitem} 
\usepackage{amsfonts,dsfont,mathtools, mathrsfs,amsthm} 
\usepackage[amssymb, thickqspace]{SIunits}
\usepackage{fancyhdr,mdframed,nicefrac}

\usepackage[norelsize, ruled, vlined]{algorithm2e}
\usepackage{algcompatible}

\usepackage{epsfig}
\usepackage{graphicx}
\usepackage{float}
\usepackage{caption}
\usepackage{subcaption}

\usepackage{courier}

\usepackage{multirow}
\usepackage{bigstrut}

\allowdisplaybreaks
\date{\today}
\makeatletter
\def\@settitle{\begin{center}%
		\baselineskip14\p@\relax
		\normalfont\LARGE\scshape\bfseries
		\@title
	\end{center}%
}

\def\@setauthors{%
  \begingroup
  \def\thanks{\protect\thanks@warning}%
  \trivlist
  \centering\footnotesize \@topsep30\p@\relax
  \advance\@topsep by -\baselineskip
  \item\relax
  \author@andify\authors
  \def\\{\protect\linebreak}%
  \authors%
  \ifx\@empty\contribs
  \else
    ,\penalty-3 \space \@setcontribs
    \@closetoccontribs
  \fi
  \endtrivlist
  \endgroup
}

\makeatother
\makeatletter

\def\subsection{\@startsection{subsection}{2}%
	\z@{.5\linespacing\@plus.7\linespacing}{.5\linespacing}%
	{\normalfont\large\bfseries}}

\def\subsubsection{\@startsection{subsubsection}{3}%
	\z@{.5\linespacing\@plus.7\linespacing}{.5\linespacing}%
	{\normalfont\itshape}}

\usepackage{multirow}
\usepackage{bigstrut}
\usepackage{wrapfig}
\usepackage{caption}

\newtheorem{theorem}{Theorem}[section]

\renewcommand{\geq}{\geqslant}
\renewcommand{\ge}{\geqslant}
\renewcommand{\le}{\leqslant}
\renewcommand{\leq}{\leqslant}

\renewcommand{\mapsto}{\longmapsto}

\DeclareSymbolFont{symbolsC}{U}{pxsyc}{m}{n}

\usepackage{bm}
\usepackage{amsthm}

\newcommand{\h}{\mathcal{H}}
\newcommand{\inv}{^{-1}}
\newcommand{\mom}{\mathrm{m}}
\usepackage[utf8]{inputenc} 
\usepackage[T1]{fontenc}    
\usepackage{hyperref}       
\usepackage{url}            
\usepackage{booktabs}       
\usepackage{amsfonts}       
\usepackage{nicefrac}       
\usepackage{microtype}      
\usepackage{amsmath}

\usepackage{mathrsfs}  
\usepackage{bm}  
\usepackage{comment}

\title{Fisher entropic Fokker--Planck model of monatomic rarefied gases}
\author{Veronica Montanaro, Lukas Netterdon, Manuel Torrilhon and Hossein Gorji}
\thanks{Corresponding author: Veronica Montanaro. \\
Emails: \href{mailto:veronica.montanaro@epfl.ch}{\texttt{veronica.montanaro@epfl.ch}}, 
\href{mailto:netterdon@acom.rwth-aachen.de}{\texttt{netterdon@acom.rwth-aachen.de}}, 
\href{mailto:mt@acom.rwth-aachen.de}{\texttt{mt@acom.rwth-aachen.de}},
\href{mailto:Mohammadhossein.Gorji@empa.ch}{\texttt{Mohammadhossein.Gorji@empa.ch}}.
Veronica Montanaro: CSQI Chair, École Polytechnique Fédérale de Lausanne, Switzerland and Laboratory for Computational Engineering, Swiss Federal Laboratories for Materials Science and Technology(EMPA), CH-8600 D\"ubendorf, Switzerland. Lukas Netterdon and Manuel Torrilhon: Chair of Applied and Computational Mathematics, RWTH Aachen University,  D-52062 Aachen, Germany. Hossein Gorji: Laboratory for Computational Engineering, Swiss Federal Laboratories for Materials Science and Technology(EMPA), CH-8600 D\"ubendorf, Switzerland}

\newcommand{\dd}[2]{ \cfrac{ \partial #1 }{ \partial #2 }}
\newcommand{\ddc}[2]{ \frac{ \partial #1 }{ \partial #2 }}

\begin{document}
\begin{abstract}
    \noindent Particle-based stochastic approximations of the Boltzmann equation are popular tools for simulations of non-equilibrium gas flows, 
    for which the Navier\textendash Stokes\textendash Fourier equations fail to provide accurate description. However, these numerical methods are computationally demanding, especially in the near-continuum regime, where the collisions become overwhelming. On the other hand, the Fokker\textendash Planck kinetic models offer an efficient alternative, as the binary collisions are described by a diffusive process. Despite the intuitive advantage, rigorous and efficient Fokker\textendash Planck approximations of the Boltzmann equation remain an open problem. 
    On one hand, the moment projection of the Fokker–Planck operator should be consistent with that of the Boltzmann operator. On the other hand, the Fokker–Planck model should be constructed in such a way that the H-theorem is satisfied. The central aim of this study is fulfilling these two categorically different constraints, i.e. moment matching and entropy dissipation, within a flexible and tractable Fokker–Planck framework. To this end, we introduce a Fisher information-based entropic constraint and demonstrate that, with a suitable polynomial expansion of the drift term, it is possible to simultaneously achieve weak moment matching while honouring the H-theorem.    
    We support our theoretical result by numerical experiments on the shock problem, validating our Fisher Entropic Fokker--Planck framework. 
\end{abstract}
\maketitle
\section{Introduction}
\subsection{Background}
\noindent Rarefied gas flows appear in a large variety of applications in physics, including atmospheric reentry flights \cite{SCHOULER2020100638}, microchannel flows \cite{SHEN2003512}, and extreme ultraviolet lithography \cite{wangNumericalInvestigationFlow2024}. Their simulation is a relevant and active research area, as conventional fluid dynamics simulations based on the Navier--Stokes--Fourier (NSF) equations fail to give an accurate description of such processes \cite{knudsenGesetzeMolekularstromungUnd2006, 585795,ewartMassFlowRate2007}. The reason lies in the fact that rarefied gases are often in strong non-equilibrium conditions. This translates to extreme gradients and jumps in the thermodynamic quantities. Therefore, the closure implied by the Linear Response Theory, upon which NSF equations are built, is violated \cite{schaafFlowRarefiedGases2015, pengStatisticalFluidMechanics2023}. \newpage
\noindent There exists, however, a more accurate theoretical description for non-equilibrium gases which goes beyond the continuum \cite{chapmanMathematicalTheoryNonuniform1990} and takes into account the molecular description of the gas. Instead of following every molecule in the system, as carried out in Molecular Dynamics \cite{alderPhaseTransitionHard1957}, the main idea is to follow the statistical description of the gas kinetics by focusing on the molecular velocity distribution $f(v,x,t)$. The latter is linked to the probability of finding a particle at time $t \in {\mathbb{R}}^{+}$ and at a given position $x \in {{\mathbb{R}}^3}$ with velocity $v \in {{\mathbb{R}}^3}$. In our case we will consider the evolution of a monatomic dilute gas, for which the molecular chaos assumption (also called \textit{Stosszahlansatz}) holds. For this setting, an accurate description of the evolution of a gas in terms of the evolution of $f$ is given by the Boltzmann equation, which takes into account binary interactions between particles through an integral collision operator.
\noindent Due to the high-dimensional support of the velocity distribution function $f$ as well as the complexity of the Boltzmann collision integral, simulations based on the Boltzmann equation are prohibitively computationally expensive. To deal with the high-dimensionality, a stochastic approximation based on the Monte Carlo algorithm,\linebreak called Direct Simulation Monte Carlo (DSMC), was developed by Bird in the late 1960s \cite{birdOneDimensionalSteadyFlows1994}. The algorithm is consistent with the Boltzmann equation and reproduces experimental results accurately over the whole range of rarefaction. However, since it explicitly performs binary collisions through a jump process, its computational cost scales with the collision frequency. This becomes problematic when approaching the continuum limit, where the number of collisions grows quickly, so the range of application of DSMC remains limited to a specific set of collision frequencies associated with quite large rarefactions. To tackle this problem, simplified particle collision operators have gained attention in the community, trading accuracy with computational efficiency. Diverse classes of stochastic particle methods exist in literature, based on different formulations like the BGK equation \cite{bhatnagarModelCollisionProcesses1954, yunEntropyProductionEllipsoidal2016}, as well as the Fokker\textendash Planck (FP) operators \cite{jennySolutionAlgorithmFluid2010, gorjiFokkerPlanckModel2011, heppKineticFokkerPlanck2020, mathiaudFokkerPlanckModel2017, gorjiEntropicFokkerPlanckKinetic2021}. \\ \ \\
\noindent Within these methods, the Fokker\textendash Planck (FP) model aims to achieve a runtime which is independent of the number of collisions, by approximating the collisions through a drift-diffusion stochastic process. 
This is in contrast to the BGK model which, despite its simplified form, still employs a jump process whose evaluation needs resolving the collision frequency.
The physical motivation behind the FP model has been justified from different perspectives. For example, the Hamiltonian dynamics can be projected onto a resolved
particle leading to the Generalized Langevin Equation, as performed in the Mori\textendash
Zwanzig formalism \cite{OptimalPredictionMoria}. Alternatively, in the Kramers\textendash Moyal expansion of the Boltzmann equation, the FP
equation can be derived after keeping the first and second order terms \cite{riskenFokkerPlanckEquationMethods1996}. Because of its expressive power to describe complex distributions, the FP equation has also been used for diffusion models \cite{pmlr-v202-lai23d}. \\ \ \\
Although the design of FP operators has been widely discussed and different models have been proposed, both for monatomic \cite{jennySolutionAlgorithmFluid2010, gorjiFokkerPlanckModel2011, heppKineticFokkerPlanck2020, mathiaudFokkerPlanckModel2016, gorjiEntropicFokkerPlanckKinetic2021,gorji2015variance} and polyatomic gases \cite{gorjiEfficientParticleFokker2014, basovModelingPolyatomicGases2024, mathiaudFokkerPlanckModel2017}, there exist still open problems to be addressed in order to unlock the full potential of these methods. In particular, we want to address the open problem of having a FP model which obtains correct
transport properties along with honoring the H-theorem. The theorem, which is a fundamental property of gases following the Boltzmann equation, states that the evolution of the entropy functional has always a negative sign until the equilibrium distribution is attained. Therefore satisfaction of the H-theorem imposes a requirement on the kinetic models, for their consistency with the Boltzmann dynamics.
In practice, it has also been observed that the H-theorem consistency has direct impact on the accurate recovery of shock profiles \cite{gorjiFokkerPlanckModel2011, gorjiEntropicFokkerPlanckKinetic2021}, its absence in the FP model leads to overly diffusive shock structures. \\ \ \\ 
One notable study to enforce the H-theorem within nonlinear-drift Fokker–Planck models is the Entropic Fokker-Planck (EFP) framework \cite{gorjiEntropicFokkerPlanckKinetic2021}. In this formulation, the consistency with the H-theorem is ensured by constructing an explicit coupling between diffusion and drift coefficients. However, the resulting entropy dissipation rate lacked a clear dynamical structure: although the decay of entropy was guaranteed in sign, its evolution was not governed by a well-defined functional form.
Following this study, beyond the sign of the entropy evolution, having a well-defined rate of entropy evolution turned out to be also important for shock profiles. Specifically, the decay rate proportional to the Fisher information with respect to the equilibrium distribution turned out to give high-fidelity results \cite{gorjiEntropicFokkerPlanckKinetic2021}. Motivated by these theoretical and empirical results, we present in this paper a novel framework which honors the entropy decay proportional to the Fisher information, which we will refer to as Fisher Entropic decay rate. Our model, which we will call Fisher Entropic Fokker\textendash Planck (FE-FP), has the advantages of bringing together the H-theorem constraint and correct relaxation rates of non-equilibrium moments, while maintaining the simplicity of implementation of the previous models. 
\subsection{Contributions and organisation}
\noindent The main contributions are summarized in the following.
\begin{itemize}
\item \textit{Entropy constraint via Fisher information.} \sloppy We introduce a Fisher-information based constraint on the FP model that enforces the H-theorem with a well-defined decay rate. 
\item \textit{Well-posedness of model coefficients.} \sloppy We establish well-posedness of the linear systems underlying the model coefficients subject to the combined entropy and moment-matching constraints.
\item \textit{Efficient algorithm and validation.} \sloppy We demonstrate close agreement with DSMC and clear improvements over the cubic-drift FP model in shock-structure problem.
\end{itemize}
The remainder of the paper is organised as follows.
\begin{enumerate}
\item \sloppy Section~\ref{sec:pre} introduces the notation (Sec.~\ref{subsec:notation}) and reviews existing FP models for monatomic gases, with emphasis on linear- and cubic-drift formulations (Sec.~\ref{subsec:review-fp}).
\item Section~\ref{sec:prob} states the main problem. After recalling the H-theorem (Sec.~\ref{subsec:h-theorem}), we formulate the requirements on entropy and moment consistency (Sec.~\ref{subsec:problem-statement}) and translate them into tractable constraints on drift–diffusion coefficients (Sec.~\ref{subsec:main-idea}).
\item Section \ref{sec:FEFP} develops the proposed FE–FP model. We introduce the general formulation (Sec.~\ref{subsec:general-formulation}), establish the entropy constraint (Sec. \ref{subsec:fefp-constraint}), derive the moment consistency system (Sec.~\ref{subsec:weak-consistency}), and combine them into a well-posed linear system for the drift coefficients (Sec. \ref{subsec:linear-system}). We further determine the entropy decay time scale (Sec. \ref{subsec:time-scale}) and analyze the associated SDEs (Sec.~\ref{subsec:sdes}).
\item Section~\ref{sec:proofs} collects the proofs of the main theoretical results, including the entropy constraint (Sec.~\ref{proof:entropy-constraint}), uniqueness of the linear system (Sec.~\ref{proof:invertibility}), and well-posedness of the velocity SDE (Sec.~\ref{proof:sde-global-solution}).
\item Section~\ref{sec:num} presents the numerical framework. We derive the explicit FE–FP formulation in three dimensions (Sec.~\ref{subsec:model-3d}), propose a solver for the linear system (Sec.~\ref{subsec:computation-system}), and describe moment estimation and time integration of the SDEs (Secs.~\ref{subsection:stochastic-moments}–\ref{subsection:time-integration}). The full algorithm is summarized in Sec.~\ref{subsection:final-algorithm}.
\item Section~\ref{sec:results} reports numerical experiments on the shock problem, comparing FE–FP with DSMC benchmark and the cubic-drift FP model.
\end{enumerate}
\section{Preliminaries} \label{sec:pre}
\noindent In this section we will introduce main elements of the notation and review a class of existing FP models relevant for our study. 
\subsection{Remarks on notation} \label{subsec:notation}
\noindent Throughout the paper, the subscripts with Latin letters will denote a vector component, i.e., $x_i$ is the i-th component of $x\in\mathbb{R}^n$. The $\ell^p(\mathbb{R}^n)$ norm of $x$, for $p\in\mathbb{N}/0$, is defined as $\|x\|_p\coloneqq \left(\sum_i |x_i|^p\right)^{1/p} $. For a multi-index $\alpha=(\alpha_1,...,\alpha_n)^T\in \mathbb{N}^n$, we refer to the length as $|\alpha|\coloneqq \sum_{i=1}^n \alpha_i$ and define the monomials $x^\alpha\coloneqq\Pi_{i=1}^n x_i^{\alpha_i}$. Furthermore let
\begin{eqnarray}
I^n_r \coloneqq \left\{ \alpha\in \mathbb{N}^n : |\alpha|\leq r\right\}  
\end{eqnarray}
be the set of multi-indices with order less than $r$. Therefore a polynomial of degree less than $r$ can be described as  $\sum_{\alpha\in I_r^n} c_\alpha x^\alpha$ with real coefficient $c_\alpha$. We denote by $|I_r^n|$ the cardinality of the set $I_r^n$.\\ \ \\
\noindent
We denote the velocity by $v\in \mathbb{R}^3$ and spatial position by $x\in \mathcal{D}$, where $\mathcal{D} \subseteq \mathbb{R}^3$. Moreover, we denote by $\{H_\alpha(v)\}_{\alpha \in I^3_r}$ a set of linearly independent polynomials of the velocity. \\ \ \\ 
Suppose $\hat{f}(v|x,t)$ is the density of a probability distribution over the velocity space, and let
\begin{eqnarray}
\Delta(\mathbb{R}^3)\coloneqq\left\{\hat{f}(v|x,t):\mathbb{R}^3\times \mathcal{D}\times \mathbb{R}^{\geq 0}\to \mathbb{R}^{\geq 0}\bigg| \int_{\mathbb{R}^3}\hat{f}(v|x,t)dv=1, \forall (x,t)\in \mathcal{D}\times \mathbb{R}^{\geq 0} \right\}. \nonumber \\
\end{eqnarray}
We refer to the scaled probability density $f(v,x,t)=\rho(x,t) \hat{f}(v|x,t)$ as the velocity distribution, where $\rho(x,t)=n(x,t)m$ is the gas density, with molecular mass $m$ and number density $n(x,t)$.
We denote the moments by
\begin{eqnarray}
\mathrm{m}_\alpha^{(p)}(f,x,t)\coloneqq \langle v^\alpha \| v\|_2^p, f \rangle  
\end{eqnarray}
where $\langle p ,q  \rangle = \int_{\mathbb{R}^3} pq \ dv$. We also refer to central moments as 
\begin{eqnarray}
\tilde{\mathrm{m}}_\alpha^{(p)}(f,x,t)\coloneqq \left\langle \left(v-{\mom_1^{(0)}}/{\rho} \right)^\alpha \left\| v-{\mom_1^{(0)}}/{\rho}\right\|_2^p, f \right\rangle  \ .
\end{eqnarray}
We indicate random variables by capital letters and assume that they belong to the Hilbert space $\mathscr{H}=\mathcal{L}^2(\Omega,\mathscr{A},P)$, with the sample space $\Omega$, the $\sigma-$algebra $\mathscr{A}$ on sub-sets of $\Omega$, and the law $P$ as a function from $\mathscr{A}$ to $[0,1]$. The space includes square-integrable $\mathscr{A}$-measurable functions on $\Omega$. A random process $X_t$ is a random function in $\mathscr{H}$ indexed by time $t\ge 0$. In particular, we use $M_t,X_t : \Omega \longrightarrow \mathbb{R}^3$ for the stochastic processes corresponding to velocity and spatial position, respectively. We refer to the pairs $\{(M_t^{(i)},X_t^{(i)})\}_{i\in[1,N]}$, which are initialized as $N$ i.i.d. samples and evolve according  to given stochastic processes, as particles. 
We estimate the velocity moments at a given position and time as the following. Consider $\mathcal{D}$ to be discretized into non-overlapping cells as $\cup_j \ S (x^{(j)})$, where $x^{(j)}$ is the spatial coordinate of the center of the cell $S(x^{(j)})$. We employ  the Monte Carlo numerical estimation of the conditional expectation $\mathbb{E}[ F(M_t) | X_t=x^{(j)} ]$, with $F: \mathbb{R}^3 \to \mathbb{R}^d$ a smooth function, as
\begin{eqnarray}
    \hat{\mathbb{E}}^{(N_S)}[F(M_t)|X_t=x^{(j)}] \coloneqq \cfrac{1}{N_S} \sum_{X_t^{(i)}\in S(x^{(j)})} F(M^{(i)}_t) \ ,
\end{eqnarray}
where $N_S$ is the number of particles in $S(x^{(j)})$. Note that for notational simplicity, we often drop $x^{(j)}$ in the conditional expectation. 
\subsection{Review of Fokker\textendash Planck models} \label{subsec:review-fp}
\noindent The Boltzmann equation for the evolution of $f$ in the monatomic dilute gas setting (in the absence of external force) reads
\begin{eqnarray} \label{eq:boltz}
    \dd{f}{t} + v\cdot \nabla_x f = S^{Boltz}[f],
\end{eqnarray}
where $S^{Boltz}[f]$ is the operator describing the binary collisions of the particles 
\begin{eqnarray} \label{eq:boltz-op}
    S^{Boltz}[f]= \int_{\mathbb{R}^3 \times \mathbb{S}^2} B\left( |v - u |,\cos\hat{\omega}\right) \left( f(v^*) f(u^*) - f(v)f(u) \right)\ du\ d\hat{\sigma}   .
\end{eqnarray}
Here $(u,v)$ and $(u^*,v^*)$ are the velocities of the colliding pairs pre and post collisions, $B: \mathbb{R}^3 \times [-1,1] \longrightarrow \mathbb{R}^+$ is the collision kernel, $\hat{\omega}$ is the deflection angle, and $d\hat{\sigma}$ is the infinitesimal solid angle \cite{cercignaniChapterBoltzmannEquation2002}. For notation simplicity, we omitted the dependency on $x$ and $t$. To conclude, it is important to mention that the moments of $f$ give information on the macroscopic quantities of the gas. We can compute the density as $\mathrm{m}^{(0)}_0$, the mean velocity as $U_i = \mom_i^{(0)}/\rho$ and the temperature as $T = \ \tilde{\mom}^{(2)}_0/ ( 3n k_B )$. In the latter, we have $k_B$ as the Boltzmann constant. For notational convenience, we refer to the quantity $\theta=k_BT/m = \tilde{\mom}^{(2)}_0 / (3\rho
)$ as well. The other relevant physical quantities depending on the moments are the scalar pressure $p = \tilde{\mom}^{(2)}_0/3$, the pressure tensor $\Pi_{ij} = \tilde{\mom}^{(0)}_{ij}$, the stress tensor $\sigma_{ij} = \Pi_{ij} - 1/3p\delta_{ij}$ and the heat-fluxes $\ q_i = \tilde{\mathrm{m}}^{(2)}_i/2$. Finally, we refer to the viscosity $\mu = \mu_0 \left( T /T_0 \right)^\omega$, where $\omega$ depends on the underlying molecular model, $T_0=273$~K and $\mu_0$ the viscosity evaluated at $T_0$. \\ \ \\
\noindent
Approximating the Boltzmann operator with a simplified diffusion operator, leading to the FP equation, dates back to the works of Kirkwood \cite{kirkwoodStatisticalMechanicalTheory1946}, Lebowitz et al. \cite{lebowitzNonequilibriumDistributionFunctions1960}, and Pawula \cite{pawulaApproximationLinearBoltzmann1967}. Based on the theory of Brownian motion, these studies presented a heuristic approach to obtain a diffusive approximation of the Boltzmann collision operator. In general such approximations can  be cast in an operator of the form
\begin{eqnarray} \label{eq:fp}
    S^{FP}[f] = - \nabla_v \cdot ( A[f](t,v) \ f) + \Delta_{v,v} (D[f](t,v) \ f) .
\end{eqnarray}
where $A[f](t,v): \Delta(\mathbb{R}^3)\times [0,\infty)\times \mathbb{R}^3  \to \mathbb{R}^3$ and $D[f](t,v): \Delta(\mathbb{R}^3)\times [0,\infty) \times \mathbb{R}^3  \to \mathbb{R}^3 \otimes \mathbb{R}^3$. We often drop the explicit dependency on $t$ and $v$ and refer to $A[f]$ as the drift functional and $D[f]$ as the diffusion functional.

\subsubsection{Linear-drift FP} \label{subsec:linear-fp}
\noindent \sloppy The original diffusive approximation of the linear Boltzmann equation from \cite{kirkwoodStatisticalMechanicalTheory1946} and \cite{pawulaApproximationLinearBoltzmann1967} led to a formulation with a drift which is linear with respect to the (molecular) velocity and a diffusion proportional to the temperature. In the context of simulation of monatomic rarefied gases, Jenny et al. \cite{jennySolutionAlgorithmFluid2010} were the first ones to adopt this formulation. The resulting operator takes the form of 

\begin{eqnarray} 
    S^{FP}[f] =  \nabla_v \cdot \left( {v^\prime}/{\tau}\ f \right) + \Delta_{v,v} \left( {\theta}/{\tau}\ f\right), \label{eq:fp-linear}
\end{eqnarray}
where
\begin{eqnarray}
    A[f] = - {v^\prime}/{\tau} \label{eq:linear-drift}
\end{eqnarray}
and
\begin{eqnarray}
    D[f] = {\theta}/{\tau} \label{eq:linear-diff}\ ,
\end{eqnarray}
where $v^\prime= v - U$ is the fluctuating velocity. Here, $\tau = 2\mu/p$ is a time scale. This model formulation leads to the Ornstein-Uhlenbeck process of McKean-Vlasov-type (since $U$ and $T$ depend on $f$). The linear-drift FP satisfies conservation of mass, momentum and energy, and recovers the correct viscosity in the hydrodynamic limit. This implies correct relaxation rate for the stress tensor $\sigma$. However, the model has a major shortcoming: it fails to recover the correct relaxation rate for the heat-fluxes $q$. This leads to a wrong Prandtl number and thus a wrong heat conductivity in the hydrodynamic limit \cite{jennySolutionAlgorithmFluid2010}.
\subsubsection{Cubic-drift FP} \label{subsec:nonlinear-fp}
\noindent In the following years, FP methods with correct heat-flux relaxation have been proposed for monatomic \cite{gorjiFokkerPlanckModel2011, mathiaudFokkerPlanckModel2016} and polyatomic gases \cite{mathiaudFokkerPlanckModel2017}. In order to fulfill consistent relaxation of moments beyond the stress tensor, a fully non-linear model of the drift was postulated in \cite{gorjiFokkerPlanckModel2011}. In particular, the diffusion 
\begin{eqnarray} \label{eq:diffusion-nonlinear}
    D[f] = {\theta}/{\tau}
\end{eqnarray}
maintains the form of the linear-drift model, and $A[f]$ now has a cubic form given by
\begin{eqnarray} \label{eq:drift-cubic}
    A_i[f](v^\prime) = - \cfrac{1}{\tau} v^\prime_i + \sum_{j=1}^3 c_{ij}[f]\ v^\prime_j + \gamma_i[f] \left( \sum_{j=1}^3 v_j^\prime v_j^\prime - 3\theta \right) + \Lambda[f] \left( v_i^\prime  \sum_{j=1}^3 v_j^\prime v^\prime_j - \cfrac{2q_i}{\rho} \right), \nonumber \\
\end{eqnarray}
where $\tau$ has the same form as for the linear model, $c[f] \in \mathbb{R}^{3\times3}$ is a symmetric tensor, $\gamma[f] \in \mathbb{R}^3$ is a vector, and $\Lambda[f]$ is a scalar. While $\Lambda [f]$ was introduced as a stability term, $c [f] $ and $\gamma [f] $ are updated at every timestep through a linear system of equations. The system is constructed using weak consistency with the Boltzmann operator as a constraint, i.e.
\begin{eqnarray} \label{eq:weak-consistency}
    \langle H_\alpha(v^\prime), S^{FP}[f]\rangle\overset{!}{=} \langle H_\alpha(v^\prime), S^{Boltz}[f] \rangle, & \forall \ \alpha \in I_\alpha^3.
\end{eqnarray}
In the case of the cubic-drift \eqref{eq:drift-cubic}, \sloppy $\{H_\alpha(v^\prime)\}_{\alpha \in I_\alpha^3} = \{1, v^\prime, \|v^\prime\|^2_2 ,v^\prime \otimes v^\prime - 1/3\mathrm{tr}(v^\prime),  v^\prime \| v^\prime \|^2_2 \}$. The objective is to enforce the projection of the FP collision operator on a given set of polynomials to be similar to that resulting from the Boltzmann operator $\forall \ t \in \mathbb{R}^{\geq0}$. For this reason, the weak consistency is also known in the community as``moment matching''. This enforces both conservation laws through the moments of $1, v^\prime, \|v^\prime\|^2_2 $ and, in addition to linear-drift FP, the relaxation rates of higher-order moments (up to the heat-fluxes). This model is the state-of-the-art for FP model thanks to the simplicity of its formulation, as well as flexibility on the moments matched. We have to notice that despite the simplicity of its implementation and its moment consistency property, the cubic-drift model has no guarantee on the structure  of the entropy dynamics and its relaxation rate towards the equilibrium. As a consequence, the model can be prone to diffusive shock profiles in numerical experiments \cite{gorjiFokkerPlanckModel2011, gorjiEntropicFokkerPlanckKinetic2021}. 
\subsubsection{FP models with alternative structures} 
\label{subsec:other-contributions}
\noindent It is worthwhile to mention that there exists two alternative structures of FP models which can address the issue of entropy relaxation encountered by the cubic-drift FP model. The first one, introduced by Mathiaud and Mieussens \cite{mathiaudFokkerPlanckModel2016}, leverages the anisotropic diffusion coefficient inspired by the elliptic BGK model \cite{yunEntropyProductionEllipsoidal2016}. Despite the entropy consistency, the moment consistency could not be enforced globally as the diffusion coefficient might become negative for specific non-equilibrium conditions. Another result was derived in \cite{gorjiEntropicFokkerPlanckKinetic2021}, where a scalar diffusion was constructed ensuring consistency with the Boltzmann H-theorem. The improvement over the cubic-drift model, however, remained limited due to the fact that the entropy relaxation rate could not be controlled.   

\section{Problem statement and main idea} \label{sec:prob}
\subsection{H-theorem} \label{subsec:h-theorem} 
\noindent Considering the Boltzmann equation in a homogeneous setting \eqref{eq:boltz} without the gradient term on the left-hand side, it can be shown that the equilibrium  is given by the Maxwell\textendash Boltzmann distribution
\begin{eqnarray} \label{eq:equilibrium}
    f_0 = \cfrac{\rho}{(2\pi \theta)^{3/2}} \exp{\left(- \cfrac{ \| v^\prime\|^2_2}{2\theta} \right)}  \ .
\end{eqnarray}
 The H-theorem characterizes the evolution of $f$ towards $f_0$ via the entropy functional
 \begin{eqnarray} \label{eq:entropy}
    \h[f] = \langle f,\log f\rangle,
\end{eqnarray}
 stating that
\begin{eqnarray} \label{eq:h-theorem}
    \dd{\h[f]}{t} \leq 0, 
\end{eqnarray}
with equality $\partial \mathcal{H}[f]/\partial t=0$ if and only if $f=f_0$. This is an important property of the Boltzmann equation, whose physical interpretation gives rise to irreversibility of the kinetic system, giving an ``arrow of time" towards the equilibrium, which corresponds to the minimum entropy $ \h [f_0]=  \langle f_0 ,\log f_0\rangle $ \cite{toscaniHtheoremAsymptoticTrend1987}. Moreover, existence of $L_1$ and Kullback-Leibler bounds with respect to the equilibrium can be proved based on contraction of $\mathcal{H}[f]$ \cite{toscaniSharpEntropyDissipation1999}. 

\subsection{Problem statement} \label{subsec:problem-statement}
\noindent Given the importance of the H-theorem in Boltzmann kinetics, it is highly desirable to devise FP models admitting this constraint. In particular, theoretical results on the trend to equilibrium in terms of entropy bounds of \eqref{eq:fp} have been shown in \cite{arnoldCONVEXSOBOLEVINEQUALITIES2001,arnoldLargetimeBehaviorNonsymmetric2008, toscaniEntropyProductionRate1999}. In practice, admitting the H-theorem has important impact on the ability of the model to accurately recover the shock profiles \cite{gorjiFokkerPlanckModel2011, gorjiEntropicFokkerPlanckKinetic2021}. \\ \ \\
In the following, we translate the H-theorem to a model constraint for a generic FP operator. Considering the homogeneous setting, the evolution of $\h[f]$ following the FP equation \eqref{eq:fp} and assuming $D[f]$ to be independent of $v$, we get
\begin{eqnarray} \label{eq:entropycalc}
\dd{\h[f]} t&=&\langle  -\nabla_v A[f] - D[f] \ (\nabla_v \log f\cdot \nabla_v \log f) , f \rangle \ .
\end{eqnarray}
\noindent
The open question is how to choose $A[f]$ and $D[f]$ such that $\partial_t \h[f]$ follows the constraints given by the Boltzmann H-theorem, and at the same time we keep enough flexibility to fulfill the weak consistency relationship \eqref{eq:weak-consistency} for an arbitrary set of moments. 
\subsection{Main idea} \label{subsec:main-idea}
\noindent Considering the linear-drift FP model, a key insight can be gained. 
By using $A[f]$ \eqref{eq:linear-drift}, $D[f]$ \eqref{eq:linear-diff} from the operator \eqref{eq:fp-linear} into the expression \eqref{eq:entropycalc} we find
\begin{eqnarray} \label{eq:linear-entropy}
 \cfrac{\partial \h[f]}{\partial t}&=&- \cfrac{\theta}{\tau} \left\langle \nabla_v \log \cfrac{f}{f_0} \cdot \nabla_v \log \cfrac{f}{f_0}, f \right \rangle  
 \label{calc:linear:end} ,
\end{eqnarray}
where the integral term on the right-hand-side is the Fisher information $\mathcal{I}(f|f_0)=\left\langle \nabla_v \log {f}/{f_0} \cdot \nabla_v \log  {f}/{f_0}, f \right\rangle$ \cite{fisherTheoryStatisticalEstimation1925}.  
The quantity $\mathcal{I}(f|f_0)$ represents a distance between two probability distributions $f$ and $f_0$. Because of its quadratic form, it is always $\geq0$ with equality iff $f=f_0$. Since $T>0$, then \eqref{eq:h-theorem} holds and the H-theorem is verified. This formulation is particularly favorable, as it is consistent with the theoretical analysis of FP equations, as seen in \cite{toscaniEntropyProductionRate1999}. Moreover it was observed that fulfilling the equality \eqref{eq:linear-entropy} leads to better recovery of shock profiles with respect to other rates \cite{gorjiEntropicFokkerPlanckKinetic2021}. However, we remind that despite this nice property, the linear-drift FP model does not give correct heat-conductivity, as shown in Sec.~\ref{subsec:linear-fp}. Therefore the main idea of our Fisher Entropic FP model is to bring together the entropy evolution \eqref{eq:linear-entropy} and the weak consistency relation \eqref{eq:weak-consistency} in the construct  of drift and diffusion coefficients of the FP model. 
\section{Entropy guided Fokker\textendash Planck by Fisher information}  \label{sec:FEFP}
\noindent This section is dedicated to discuss in detail the FE-FP framework. Closures of the diffusion and drift functionals, as well as important theoretical results on the model, will be covered. 
The general ansatz for the drift will be introduced in Sec.~\ref{subsec:general-formulation}. The novel entropy constraint will be discussed in Sec.~\ref{subsec:fefp-constraint}. The moment consistency relations will be then derived in Sec.~\ref{subsec:weak-consistency}. We will then proceed to Sec.~\ref{subsec:linear-system} to show how to include both constraints into a linear system for which there exists a unique solution. Finally, we will go through the time scale derivation in Sec.~\ref{subsec:time-scale} and conclude by analyzing the underlying SDE in Sec.~\ref{subsec:sdes}.
\subsection{General formulation} \label{subsec:general-formulation}
\noindent We start sketching our model by constructing $A$ as a gradient field which, without loss of generality, can be cast as an expansion around the linear drift. 
Suppose we would like to match the moments associated with polynomials $\{H_\alpha(v^\prime)\}_{\alpha \in I_r^3}$. Consider $\{H_\alpha(v^\prime)\}_{\alpha \in I_r^3} \cup \Vert v^\prime \Vert_2^{2k}$, where $2k=r+2-\textrm{mod}(r,2)$. The augmented set of polynomials now include the smallest even power, larger than the highest order of moments which are going to be matched. We refer to the augmented index set by $\hat{I}_{2k}^3$ and the corresponding polynomials by $\{H_\alpha(v^\prime)\}_{\alpha \in \hat{I}_{2k}^3}$.
Let $\Phi^{(r)}[f](v^\prime)$ be a potential function corresponding to the polynomial expansion, given by
\begin{eqnarray} \label{eq:phi}
     \Phi^{(r)}[f](v^\prime) \coloneqq \sum_{\alpha \in \hat{I}^3_{2k}} c_\alpha[f] \  H_\alpha(v^\prime)
\end{eqnarray}
with $\{c_\alpha[f]\}_{\alpha \in \hat{I}^3_{2k}}$ being a set of real coefficients.
We consider the following drift ansatz
\begin{eqnarray} \label{eq:fe-fp-drift}
    A[f](v^\prime) = -\cfrac{1}{\tau}v' + A^{HO}[f](v^\prime) - c_d[f] \nabla_v (\|v^\prime\|^d_2), 
\end{eqnarray}
comprising of linear, nonlinear, and stabilizing term, where
\begin{eqnarray}
     A^{HO}[f](v^\prime) = \nabla_v \Phi^{(r)}[f](v^\prime) \ ,
\end{eqnarray} 
 contains all higher order terms in $v^\prime$, up to the order $r$. Note that $\tau>0$ is a time scale independent of $v$ that will be fixed later.
 The diffusion is kept with the same formulation as the one expressed in \eqref{eq:diffusion-nonlinear}. \\ \ \\
 \noindent
The purpose of the stabilizing term $-c_d[f]\ \nabla_v \Vert v^\prime \Vert_2^d $ is to ensure that the drift can be expressed as a gradient of a convex function as $\Vert v^\prime \Vert_2 \to \infty$. This leads to regularity of the solution of the corresponding SDE (as justified in Sec.~\ref{subsec:sdes}). We choose $c_d[f] \geq 0$ (with equality at equilibrium),
where $d = 2k + 2$ is the closest even number strictly larger than $2k$. Furthermore, as we would like to recover the linear drift model in the limit of thermal equilibrium, the coefficient $c_d[f]$ should depend on some notion of non-equilibrium moments. We set $c_d[f]$  in terms of $\theta_0$ and a small empirical parameter $\epsilon_0 \approx 10^{-3}$ as
\begin{eqnarray} 
c_d[f] &=& \epsilon_0 \sqrt{(2\theta_0)^r}  \ 
     \left\| \tilde{\mathrm{m}}^{(r)}_0(f)\ - \tilde{\mathrm{m}}^{(r)}_0(f_0)  \right\|^2_2 \ \Big/\ \left\|  \tilde{\mathrm{m}}^{(r)}_0(f_0)  \right\|^2_2 .
    \label{def:stab-coeff}
\end{eqnarray} 
 The idea, following the one for the cubic-drift formulation of Sec.~\ref{subsec:nonlinear-fp}, is to design a linear system of equations for $\{c_\alpha[f]\}_{\alpha \in \hat{I}_{2k}^3}$ such that moment consistency and entropy production are both satisfied. 
\subsection{Fisher entropic constraint} \label{subsec:fefp-constraint}
\noindent We want to derive a constraint which recovers the equality \eqref{eq:linear-entropy} given the formulation of the drift in \eqref{eq:fe-fp-drift}. By carrying out the calculation as in \eqref{eq:entropycalc}, we arrive to the following result (the proof is deferred to Sec.~\ref{proof:entropy-constraint}).
\begin{theorem}[Fisher entropic constraint] \label{theorem:entropy-constraint}
    Suppose $f(v,t)$ (with $\mathcal{H}[f]<\infty$) is the solution of the FP equation 
    \begin{eqnarray}
    \frac{\partial}{\partial t} f&=&S^{FP}[f] \ ,
    \end{eqnarray}
    for some fixed initial condition, where the right-hand-side is given by \eqref{eq:fp}, the drift follows \eqref{eq:fe-fp-drift}, and the diffusion is set to \eqref{eq:diffusion-nonlinear}. Also for simplicity (and without loss of generality) $\rho=1$. If the drift coefficients satisfy

    \begin{eqnarray}&& \label{eq:entropy-constraint}
    \left\langle \sum_{\alpha \in \hat{I}_{2k}^3} c_\alpha[f]\ \Delta_{v,v} H_\alpha(v^\prime), f\right\rangle = (2k+2)(2k+1)c_d[f] \langle \|v^\prime\|^{2k}_2, f\rangle
    \end{eqnarray}
\noindent then, the following holds    
    \begin{eqnarray} \label{eq:fisher-entropy}
            \dd{\h[f]}{t} =  -D[f]\ \mathcal{I}(f \ |\ f_0).
    \end{eqnarray}
\noindent Here, $\mathcal{I}$ is the Fisher information defined as
    \begin{eqnarray}
        \mathcal{I}(f\ |\ f_0) \coloneqq \left\langle \nabla_v \log \cfrac{f}{f_0} \cdot  \nabla_v \log \cfrac{f}{f_0}, f \right\rangle \ .
    \end{eqnarray} 
\end{theorem}
\noindent
Theorem \ref{theorem:entropy-constraint} rests at the core of the devised model, since it shows
that by enforcing \eqref{eq:entropy-constraint} on the drift coefficients, we obtain the
favorable entropy decay given by \eqref{eq:fisher-entropy}. Note that the time scale $\tau$, introduced
in $D[f]$, controls the entropy relaxation rate. 
While well-posedness and convergence to equilibrium of the FE-FP system remains out of the scope of the current study, we emphasize that the Fisher entropic structure imposed by \eqref{eq:fisher-entropy} leads to an upper bound on the $L^1$-distance of $f$ from $f_0$ as shown for non-linear FP equations in \cite{arnoldCONVEXSOBOLEVINEQUALITIES2001}. 
\subsection{Weak consistency} \label{subsec:weak-consistency}
\noindent After deriving the entropy constraint of our FE-FP model, we move on to the moment equations resulting from the weak consistency constraints \eqref{eq:weak-consistency}. \\ \ \\ \noindent
An important result regarding relaxation of non-equilibrium moments according to the Boltzmann equation (assuming molecules with the Maxwell potential) is the theorem of the trend to equilibrium \cite{truesdellTrendEquilibriumAccording1984a}, recalled below.  Consider an arbitrary (regardless of whether it is central) moment of the form  $ \langle H_\alpha(v^\prime), f \rangle$, whose evolution is governed by the homogeneous Boltzmann equation leading to
\begin{eqnarray} \label{eq:production-term}
  \dd{}{t} \langle H_\alpha(v^\prime), f \rangle =  \left\langle  H_\alpha(v^\prime),S^{Boltz}[f]\right\rangle =: P_\alpha[f] ,
\end{eqnarray}
where $P_\alpha[f]$ is called the production term for $H_\alpha(v^\prime)$. Then we have the following result.
\begin{theorem}[trend to equilibrium] \label{theorem:prod-terms}
   Every moment following ODE \eqref{eq:production-term}  is analytic function of time which relaxes exponentially to its equilibrium value. 
\end{theorem}
\noindent For the proof, see the reference text \cite{truesdellTrendEquilibriumAccording1984a}. 
The consequence of this theorem is that the system of equations constructed from the constraint \eqref{eq:weak-consistency} is closed for any order of moments. Following this theoretical result, we insert the FP operator \eqref{eq:fp} with drift \eqref{eq:fe-fp-drift} and diffusion \eqref{eq:diffusion-nonlinear} in the left-hand side of \eqref{eq:weak-consistency}, leading to 
\begin{eqnarray} 
    \left\langle  \nabla_vH_\alpha(v^\prime) \cdot  \sum_{\beta \in I_r^3} c_\beta[f]\ \nabla_v H_\beta(v^\prime) , f \right\rangle = P_\alpha[f] + g_\alpha(\tau,D[f],c_d[f])\ , \ \forall \alpha \in I_r^3. \nonumber \\
    \label{eq:moment-equations}
\end{eqnarray}
Here, we have grouped all the right-hand side terms depending on $D[f]$, $\tau$ and $c_d[f]$ in the function 
\begin{eqnarray}
    g_\alpha(\tau,D[f],c_d[f]) &=&  \cfrac{1}{\tau}\ \langle v' \nabla_v H_{\alpha}(v^\prime), f \rangle  - D[f] \langle \Delta_{v,v}  H_{\alpha}(v^\prime) , f \rangle \nonumber \\
    &&-c_d[f] \langle \nabla_v \|v^\prime\|^d_2 \cdot \nabla_v H_\alpha(v^\prime), f\rangle.
\end{eqnarray}

\subsection{Linear system for unknown coefficients} \label{subsec:linear-system}
\noindent We have now the necessary ingredients to construct the linear system for the coefficients of the FE-FP model $\{c_\alpha[f]\}_{\alpha \in \hat{I}_{2k}^3}$, which include the moment matching constraints \eqref{eq:entropy-constraint} and the entropy constraint \eqref{eq:moment-equations}. The construction requires special caution, though, since it is not always guaranteed that the matrix obtained from combining $\Phi^{(r)}[f](v^\prime)$ with the constraints is invertible. Therefore, in order to obtain a well-posed system we proceed with the following procedure.
\begin{enumerate}
    \item We split the model coefficients $\{c_\alpha[f] \}_{\alpha \in \hat{I}_{2k}^3}$ into two sets of coefficients, and derive the potential spanned by them.
    \item We show that by enforcing the model constraints on the new coefficients, two well-posed linear systems are obtained.
\end{enumerate}
This procedure will be such that the existence of a unique solution is ensured, with a minimal additional computation of the two sets of model coefficients. \\ \ \\
Let $\varphi\coloneqq\alpha_{2k}$, $l\coloneqq |\hat{I}_{2k}^3 |$ and note that $|I_r^3 | = l-1$. Furthermore, consider the matrix $R[f] \in \mathbb{R}^{l\times l}$ with elements
\begin{eqnarray}
    R_{\alpha \beta}[f] = \langle  \nabla_v H_{\alpha}(v^\prime)\ \nabla_v H_{\beta}(v^\prime),f \rangle, \label{def:R}
\end{eqnarray}
and the vector $Q[f] \in \mathbb{R}^{l}$, 
\begin{eqnarray} \label{def:q}
    Q_\alpha[f] = \langle \Delta_{v,v} H_{\alpha}(v^\prime) , f \rangle
\end{eqnarray}
with $\alpha, \beta \in \hat{I}_{2k}^3$.
Additionally, let $G[f] \in  \mathbb{R}^l$ be a vector of zeros except for the last term
\begin{eqnarray}
    G_{\varphi}[f] = -2Q_{\varphi}[f].
\end{eqnarray}
Let us introduce $\{c^\prime_{\alpha} [f] \}_{\alpha\in \hat{I}_{2k}^3}$ and $\{\hat{c}_\alpha [f] \}_{\alpha\in \hat{I}_{2k}^3}$, from which $\{c_\alpha[f]\}_{\alpha\in \hat{I}_{2k}^3}$ are computed
\begin{eqnarray} \label{eq:coeff-relation}
    c_\alpha[f] &\coloneq& c^\prime_\alpha[f] + c^\prime_{\varphi}[f]\ \hat{c}_\alpha [f]
, \  \forall\ \alpha \in {I}_r^3 \nonumber \\
     c_{\varphi}[f] &\coloneq& c^\prime_{\varphi}[f]\ \hat{c}_{\varphi}[f] \ .
\end{eqnarray}
Accordingly, an equivalent formulation of the drift can be obtained by 
\begin{eqnarray} \label{eq:phi-proxy}
     \tilde{\Phi}^{(r)}[f](v^\prime) = \sum_{\alpha \in {I}_{r}^3} c^\prime_\alpha[f]\ H_\alpha(v^\prime) + c^\prime_{\varphi}[f] \sum_{\alpha \in \hat{I}_{2k}^3} \hat{c}_\alpha[f]\ H_\alpha(v^\prime),
 \end{eqnarray}
 leading to $\nabla_v\Phi^{(r)}[f](v^\prime) = \nabla_v\tilde{\Phi}^{(r)}[f](v^\prime)$. 
It is now left to drive, and to show the invertibility of, the linear system of equations for $\{c^\prime_{\alpha}[f]\}_{\alpha\in \hat{I}_{2k}^3}$ and $\{\hat{c}_\alpha[f]\}_{\alpha\in \hat{I}_{2k}^3}$, fulfilling moment and entropy constraints. 
We set
\begin{eqnarray} \label{eq:c-tilde}
    \hat{c}[f] \coloneqq R[f] \inv (Q[f] + G[f] ) \ ,
\end{eqnarray}
noting that $R[f]$ is invertible as shown in Appendix~\ref{appendix:a}. 
Next, we identify $\{c^\prime_\alpha[f]\}_{\alpha\in \hat{I}_{2k}^3}$, by plugging the potential $\tilde{\Phi}^{(r)}[f](v^\prime)$ in the moment and entropy constraints, and substituting $\{\hat{c}_\alpha[f]\}_{_{\alpha\in \hat{I}_{2k}^3}}$ with its solution \eqref{eq:c-tilde}. We obtain 
\begin{eqnarray} \label{eq:linear-system}
     c^\prime[f] = L[f] \inv b[f]
\end{eqnarray}
with
\begin{eqnarray} \label{def:matrix}
    L[f] = \begin{pmatrix}
        \bar{R}[f] & \bar{Q}[f] \\
        \bar{Q}[f]^T & S[f]
    \end{pmatrix} & \ \ \textrm{and} \ \ & b[f] = \begin{pmatrix}
        \overline{b}[f] \\
        h[f] 
    \end{pmatrix} \ ,
\end{eqnarray}
where the composition of $L[f]$ and $b[f]$ is given in the following. $\bar{R}[f] \in \mathbb{R}^{l-1\times l-1}$ is the sub-matrix of $R[f]$ excluding the row $\{R_{\varphi\beta}[f]\}_{\beta \in \hat{I}_{2k}^3}$ and the column $\{R_{\alpha\varphi}[f]\}_{\alpha \in \hat{I}_{2k}^3}$. In the same way, $\bar{Q}[f]  \in \mathbb{R}^{l-1}$ is the sub-vector of $Q$ excluding its last element $Q_{\varphi}[f]$. Furthermore, we have the scalar
\begin{eqnarray}
    S[f] &=& \sum_{\alpha \in \hat{I}^3_{2k}} \hat{c}_\alpha[f]\ Q_\alpha[f] = \sum_{\alpha,\beta \in \hat{I}_{2k}^3} R_{\alpha\beta}[f]\inv Q_\alpha[f] (Q_\beta[f] + G_\beta[f]) \nonumber \\
    &=& \sum_{\alpha,\beta \in {\hat{I}_{2k}^3}} R_{\alpha\beta}[f]\inv Q_\alpha[f] Q_\beta[f]\ - 2Q_{\varphi}[f] \sum_{\alpha,\beta \in \hat{I}_{2k}^3} R_{\varphi\alpha}[f]\inv Q_\beta[f].
\end{eqnarray}
The right-hand-side vector is comprised of
\begin{eqnarray}
    \overline{b}_\alpha[f] = P_{\alpha}[f]\ + g_\alpha(\tau,D[f], c_d[f]),  \ \alpha \in {I}_r^3
\end{eqnarray}
containing the weak consistency right-hand-side, and
\begin{eqnarray}
   h[f] = -c_d[f] \langle \|v^\prime\|^{2k}_2, f\rangle \ ,
\end{eqnarray}
the entropy constraint right-hand-side. \\ \ \\
With these elements set, we have the following result. 
\begin{theorem}[linear systems solution] \label{theorem:invertibility}
 The linear systems given by \eqref{eq:c-tilde} and \eqref{eq:linear-system} always admit a unique solution as long as temperature is positive i.e. $T>0$.
 \end{theorem}
 \noindent The proof is deferred to Sec.~\ref{proof:invertibility}.
 Since both linear systems for $\{c^\prime_{\alpha}[f]\}_{\alpha\in \hat{I}_{2k}^3}$ and $\{\hat{c}_\alpha[f]\}_{\alpha\in \hat{I}_{2k}^3}$ have a unique solution, then the coefficients $\{c_\alpha[f]\}_{\alpha\in \hat{I}_{2k}^3}$ are uniquely identified through \eqref{eq:coeff-relation}, resulting in a drift fulfilling both moment and entropy constraints. 
\subsection{Time scale derivation} \label{subsec:time-scale}
\noindent The model will be closed once we fix the time scale $\tau >0$.  Note that irrespective of the value of the time scale, the model fulfills the entropy and moment constraints. Our derivation of $\tau$ will be guided by the fact that it controls the entropy decay rate. In particular, we aim for a value of $\tau$ such that the entropy decay rate implied by the FE-FP model is identical to the one resulting from the Boltzmann equation for a known class of distributions. In order to do that, we limit our analysis to the case for which an analytical solution can be obtained in time. \\ \ \\
We focus on the relaxation of elliptic (anisotropic-Gaussian) distributions of the form
\begin{eqnarray}
      f_G = \cfrac{\rho}{(2\pi)^{3/2}\sqrt{\det \Pi }} \exp{\left( -\cfrac{1}{2}\sum_{i,j} v'_{i} \Pi\inv_{ij} v'_j \right)}.
\end{eqnarray} 
To simplify the derivation, and without loss of generality, we perform the calculation on a frame of reference in which $\Pi$ is diagonal with components $\lambda_1$, $\lambda_2$ and $\lambda_3=3-\lambda_1-\lambda_2$, assuming that $\theta=1$ and $\rho=1$. 
We derive the entropy decay rate for each of the operators
\begin{eqnarray} 
    \dd{\h^{Boltz}[f_G]}{t} &=& \left\langle 1, \ddc{}{t} ( f_G \log f_G) \right\rangle, \label{eq:entropy-boltz} \\
    \dd{\h^{FE-FP}[f_G]}{t} &=& - \cfrac{1}{\tau}  \left\langle \nabla_v \log\left(\cfrac{f_G}{f_0}\right)\cdot \nabla_v \log\left(\cfrac{f_G}{f_0}\right),  f_G\right\rangle. \label{eq:entropy-fefp}
\end{eqnarray}
The former, given by \eqref{eq:entropy-boltz}, follows
\begin{eqnarray}
    \left\langle 1, \ddc{}{t} ( f_G \log f_G) \right\rangle &=&  \left\langle \dd{f_G}{t}, 1+ \log f_G \right\rangle  \nonumber \\
    &=& \sum_{i=1}^3\dd{\lambda_i}{t} \left\langle \dd{f_G}{\lambda_i} , 1 + \log f_G \right\rangle \nonumber \\
    &=& - \cfrac{p}{2\mu} \sum_{i=1}^3\left( \cfrac{1}{\lambda_i} -1 \right), \label{calc:boltz-rate}
\end{eqnarray}
where we used the fact that $\partial_t \Pi =  \partial_t \sigma=-(p/\mu) \sigma$
for interaction of Maxwell molecules governed by the Boltzmann equation. On the other hand, the entropy decay rate given by \eqref{eq:entropy-fefp} follows
\begin{eqnarray}
    -\cfrac{1}{\tau} \left\langle \nabla_v \log\left(\cfrac{f_G}{f_0}\right) \nabla_v \log\left(\cfrac{f_G}{f_0}\right), f_G \right\rangle =  - \cfrac{1}{\tau} \sum_{i=1}^3\left( \cfrac{1}{\lambda_i} - 1 \right)^2 \lambda_i\ . \label{calc:fefp-rate}
\end{eqnarray}
Enforcing equality between the two entropy decay rates leads to 
\begin{eqnarray} \label{eq:timescale-setup}
    \tau = \cfrac{2\mu}{p} \sum_{i=1}^3\cfrac{\left( {1/\lambda_i} - 1 \right)^2\lambda_i}{\left( 1/{\lambda_i} - 1 \right)}=\frac{2\mu}{p} \ .
\end{eqnarray}
\subsection{Well-posedeness of corresponding SDEs} \label{subsec:sdes}
\noindent Following standard results in theory of stochastic processes, it is straight-forward to derive the Langevin SDEs corresponding to the FP equation given by \eqref{eq:fp}. Specifically, we have
\begin{eqnarray} \label{eq:langevin-eq}
    dM_{t,i} &=& A_i[f]dt\ + \sqrt{2D[f]}\ dW_{t,i}, \\
    dX_{t,i} &=& M_{t,i}dt, \label{eq:pos-langevin-eq}
\end{eqnarray}
where $W_t : \Omega \to \mathbb{R}^3$ is the Wiener process. For our FE-FP setting, the system belongs to a family of SDEs known as McKean-Vlasov (MKV) SDEs, for which the evolution of the stochastic processes depend also on the distribution $f$ (here through moments of $f$)\cite{mckeanCLASSMARKOVPROCESSES1966, chaintronPropagationChaosReview2022}. It is easy to see that in our case, this dependency is a direct consequence of the dependency on the moments for $A[f]$ and $D[f]$. This carries additional challenges compared to classical SDEs. In particular, it is more complex to prove bounds on stability and results on well-posedness. It is out of the scope of the paper to prove regularity of our SDE in the general MKV setting . Instead, we focus on regularity of the velocity SDE with frozen model coefficients. In particular, we assume the coefficients $\{c_\alpha\}_{\alpha \in \hat{I}^3_{2k}}$ as well as the diffusion $D$ to be constant in time. For this setting, we have the following theoretical result.


\begin{theorem}[SDE regularity] \label{theorem:sde-global-solution}
Let us consider the velocity SDE \eqref{eq:langevin-eq}, corresponding to the FP equation \eqref{eq:fp} with diffusion \eqref{eq:diffusion-nonlinear} and drift defined as \eqref{eq:fe-fp-drift}, with $\{c_\alpha\}_{\alpha \in \hat{I}^3_{2k}}$ bounded and constant in time, and stabilization coefficient $c_d \geq 0$ defined as \eqref{def:stab-coeff}. Then, the velocity SDE admits a unique solution global in time.    
\end{theorem}
\noindent The proof is deferred to Sec.~\ref{proof:sde-global-solution}.

\section{Justifications of theoretical results} \label{sec:proofs}
\noindent We present in this section the proofs of the theorems justifying our theoretical results on FE-FP. Only for this section, for notation simplicity and without loss of generality, we adopt the following assumptions.
\begin{enumerate}
    \item We set $\langle v,f \rangle = 0$, such that $v^\prime = v$.
    \item We set $\rho=1$.
\end{enumerate}
From this point onwards, we suppress the explicit dependence on $f$ for notational convenience, unless required for clarity.
\subsection{Proof of Fisher entropic constraint theorem (Theorem \ref{theorem:entropy-constraint})} \label{proof:entropy-constraint}
\begin{proof}
    The proof follows from the evaluation of \eqref{eq:entropycalc}. For the drift given by \eqref{eq:fe-fp-drift} we find
\begin{eqnarray} 
   \dd{\h}{t} &=& \cfrac{3}{\tau}-\langle \nabla_v A^{HO}(v), f\rangle + c_d  \langle \Delta_{v,v}\|v\|^{d}_2, f\rangle - \cfrac{\theta}{\tau} \langle \nabla_v \log f \cdot \nabla_v \log f, f \rangle \nonumber \\
    &=& -\cfrac{\theta}{\tau}\ \mathcal{I}(f|f_0)-\langle \nabla_v A^{HO} (v), f\rangle + c_d  \langle \Delta_{v,v}\|v\|^{d}_2, f\rangle   
\end{eqnarray}
using integration by part and direct calculation. 
The constraint \eqref{eq:fisher-entropy} leads to cancellation of the second and the third term on the right-hand-side, and thus
\begin{eqnarray}
\dd{\h }{t} &=& -\cfrac{\theta}{\tau} \ \mathcal{I}(f|f_0) \ .
\end{eqnarray}
\end{proof}

\subsection{Proof of linear systems solution theorem  (Theorem \ref{theorem:invertibility})} \label{proof:invertibility}
\begin{proof}
   First, note that for the system \eqref{eq:c-tilde}, the proof is straight-forward as $R $ is symmetric positive definite (see Appendix~\ref{appendix:a} for details).
Next, for the solution uniqueness of \eqref{eq:linear-system}, we need to prove that $L $ is invertible. The idea is to use the block construction of $L $ and prove invertibility through Sylvester's criterion. \\ \ \\  
Since $R$ is symmetric positive definite, and thus invertible as previously mentioned, we can compute the determinant of $L $ via the Schur complement of $\bar{R} $
\begin{eqnarray}
    L /\bar{R}  = S  - \bar{Q} ^T\bar{R} \inv \bar{Q} .
\end{eqnarray}
\noindent  Note that
\begin{eqnarray}
    \bar{Q} ^T\bar{R} \inv \bar{Q}  &=& \sum_{\alpha,\beta \in {I}_{r}^3} R_{\alpha\beta} \inv Q_\beta  Q_\alpha  \nonumber \\
    &=& \sum_{\alpha,\beta \in \hat{I}^3_{2k}} R_{\alpha\beta} \inv Q_\beta  Q_\alpha  - \sum_{\alpha \in \hat{I}_{2k}^3} R_{\varphi \beta} \inv Q_\beta  Q_{\varphi} - \sum_{\alpha \in {I}_r^3} R_{\varphi\alpha
    } \inv Q_{\varphi} Q_\alpha  \nonumber \\
    &=& \sum_{\alpha,\beta \in \hat{I}^3_{2k}} R_{\alpha
\beta} \inv Q_\beta  Q_\alpha   - 2Q_{\varphi}  \sum_{\alpha \in \hat{I}_{2k}^3} R_{\alpha \varphi} \inv Q_\alpha  \nonumber + R_{\varphi \varphi} \inv Q_{\varphi}  Q_{\varphi}  \nonumber \\
    &=& S  + R_{\varphi \varphi} \inv Q_{\varphi}  Q_{\varphi} .
\end{eqnarray}
Therefore
\begin{eqnarray}
    L /\bar{R}  = - R_{\varphi \varphi} \inv Q_{\varphi} Q_{\varphi} .
\end{eqnarray}
However the right hand side is strictly negative and therefore the system always admits a unique solution, due to the following.
\begin{enumerate}
\item Note that $Q_{\varphi} =\langle \Delta_{v,v}(\|v\|^{2k}_2),f \rangle$ and therefore $Q_{\varphi} >0$ since $\langle \Vert v \Vert_2^2, f\rangle >0$ (due to positivity of temperature). 
\item We have $R_{\varphi\varphi} \inv>0$ by observing
\begin{equation}
R_{\varphi\varphi} \inv = \langle \nabla_v (\| v \|^{2k}_2) \cdot \nabla_v (\| v \|^{2k}_2)  , f \rangle \ .
\end{equation}
\end{enumerate}
\end{proof}
\subsection{Proof of SDE regularity (Theorem \ref{theorem:sde-global-solution})} \label{proof:sde-global-solution}
\begin{proof}
The proof follows Theorem 3.5 by Khasminskii in \cite{khasminskiiStochasticStabilityDifferential2012}. Given the fact that for fixed model coefficients, the drift is only function of $v$ and the diffusion is a constant, the conditions of the mentioned theorem on drift and diffusion are fulfilled. We are left with identifying a Lyupanov function $V(v)$ which behaves like a convex function as $\Vert v \Vert_2 \to \infty$, and satisfies
\begin{eqnarray}
L V(v) \le a V(v) \ ,
\end{eqnarray} 
where $L$ is the generator of the SDE and $a<\infty$ is some constant.
In our case, we set $V(v) = \Vert v \Vert_2^2$. We observe that as $\Vert v\Vert \longrightarrow \infty$, the stability term \eqref{def:stab-coeff} becomes dominant, as $\| v\|^d_2$ has a higher order than any other polynomial in $A^{HO} (v)$, and goes to $-\infty$ faster than any other term leading to 
\begin{eqnarray} \label{eq:generator}
         -A (v)\cdot  \nabla_{v} V(v)  + \sqrt{2\theta/\tau}\  \Delta_{v,v} V(v) \leq aV(v), 
    \end{eqnarray}
for some constant $a>0$. 
\end{proof}

\section{Solution algorithm} \label{sec:num}
\noindent For practical rarefied gas flow simulations, we need to deal with the three-dimensional velocity space, where the moment consistency includes the stress tensor and the heat-fluxes. This leads to correct viscosity and heat-conductivity in the hydrodynamic limit. In this section, we present a computationally affordable numerical recipe for the FE-FP operator built on this setting.  
\subsection{Model expression} \label{subsec:model-3d}
\noindent \sloppy In order to consider the moment consistency for the momentum, the stress tensor, and the heat-fluxes, we let $\{H_\alpha(v^\prime)\}_{\alpha \in \hat{I}^3_{4}} = \{ 1, v^\prime, v'\otimes v' - 1/3 \mathrm{tr}(v^\prime), v^\prime \|v^\prime\|^2_2, \|v^\prime\|^4_2 \}$ (note that the energy conservation is included through the trace of the stress tensor). With this set of polynomials we have the following expression for the drift
\begin{eqnarray} \label{eq:A-3d}
A_i (v^\prime) &=& -\cfrac{1}{\tau}v^\prime_i + A_i^{HO} (v^\prime) - c^{(4)} \ v^\prime_i\| v^\prime\|_2^4, 
\end{eqnarray}
where 
\begin{eqnarray}
A^{HO}_i (v^\prime)&=&c^{(0)}_i +\sum_{j=1}^3 \left( (c^{(1)}_{ij} +c^{(1)}_{ji} )v_j^\prime+2c^{(2)}_{j} \ v^\prime_{i}v^\prime_{j}\right) \nonumber \\ &+&c^{(2)}_i \ \|v^\prime\|_2^2 +3c^{(3)}  \ v^\prime_i\|v\|_2^2 .
\end{eqnarray}
For convenience we define the symmetric tensor $\tilde{c}^{(1)}_{ij} =c_{ij}^{(1)} +c_{ji}^{(1)} $.
Note that $c^{(4)}  \geq 0, |c^{(4)} | \ll 1$ is fixed as the stabilization term (see \eqref{def:stab-coeff}).
The remaining coefficients are $\{c^{(0)} , \tilde{c}^{(1)} , c^{(2)} , c^{(3)} \}$, which depend on $f$ through its moments.
The vector $c^{(0)}  \in \mathbb{R}^3$ can be readily computed from the momentum conservation 
\begin{eqnarray} 
    c^{(0)}_i  = -3\theta c^{(2)}_i  - \sum_{j=1}^3 \left( 2{\Pi_{ij}} c_j^{(2)}  + 6\cfrac{q_{i}}{\rho}c^{(3)}   \right)+\frac{\Delta_i}{\rho}c^{(4)} ,
\end{eqnarray}
with $\Delta_i  = \tilde{\mathrm{m}}^{(4)}_{i}$. The linear system is then constructed  for the symmetric tensor $\tilde{c}^{(1)}  \in \mathbb{R}^{3 \times 3}$, the vector $c^{(2)}  \in \mathbb{R}^3$, and the scalar  $c^{(3)} \in\mathbb{R}$ leading to $10\times 10$ system. 
\subsection{Linear system setup} \label{subsec:computation-system}
\noindent Following Sec.~\ref{subsec:linear-system}, we obtain the unknowns by solving the two linear systems \eqref{eq:c-tilde} and \eqref{eq:linear-system} of size $10\times10$. The explicit expressions of the system matrices $R $, $Q $, $G $ and $L $ (using the naming introduced in Sec.~\ref{subsec:linear-system}), as well as the one for the right-hand side vector, can be found in Appendix~\ref{appendix:c}.
To reduce the cost of inversion of the matrix $L $, we calculate $L \inv$ through its Schur complement (which we have shown in Sec.~\ref{proof:invertibility} to be the scalar $L/\bar{R} = \bar{R}_{10,10}\inv Q_{10}Q_{10}$)
\begin{eqnarray}
    L\inv = \begin{pmatrix}
        \bar{R}\inv-(L/\bar{R})\inv\bar{R}\inv Q^TQ\bar{R}\inv & -\bar{R}\inv \bar{Q}^T (L/\bar{R} )\inv \\
        (L/\bar{R} )\inv Q\bar{R}\inv & (L/\bar{R})\inv 
    \end{pmatrix},
\end{eqnarray}
and we get $\bar{R} \inv$ by applying the Sherman-Morrison formula \cite{shermanAdjustmentInverseMatrix1950} to $R \inv$.
In this way, we reduce the computation of the linear system to only one $10\times10$ matrix inversion, substituting the second inversion with less costly operations.
\subsection{Moment estimation} \label{subsection:stochastic-moments}
\noindent Since we pursue a particle based approach, the moments of $f$ are estimated based on the particle samples. In particular, we are interested in moments throughout  the physical space, and hence we make use of the Monte Carlo estimator for conditional expectations
\begin{eqnarray} \label{def:montecarlo-integration}
   \mathrm{m}^{(p)}_{\alpha}(f,x,t) \approx \rho(x,t) \  \hat{\mathbb{E}}^{(N)}\left[  M_t^\alpha \Vert  M_t \Vert^p_2 \Big| X_t=x\right],
\end{eqnarray}
where $N$ is the number of particles residing in the computational cell around $x$. Moreover we estimate the central moments $\mathrm{m}^{(p)}_\alpha$ with the corresponding Monte Carlo estimator for conditional expectations of $M_n^{\prime(i)} = M_n^{(i)} - \hat{\mathbb{E}}^{(N)}[M_n|X_n =x]$.
\subsection{Time integration} \label{subsection:time-integration}
\noindent We discretize the McKean-Vlasov processes \eqref{eq:langevin-eq} and \eqref{eq:pos-langevin-eq} in time, resulting in the approximation of the system of particles at $\{t_n\}$ for a given time-step size $\Delta t = t_{n+1}-t_n$. Let us denote the approximations at time $t_n$ by subscript $n$ and let $s_n=\Delta t/\tau_n$. Moreover, let $\{\xi^{(i)}\}$ be i.i.d. standard normal variates. To evolve $(M_n^{(i)},X_n^{(i)}) \mapsto (M_{n+1}^{(i)},X_{n+1}^{(i)})$, we employ the split-step time integration scheme from \cite{gorjiFokkerPlanckModel2011}, with the half-step
\begin{eqnarray} \label{eq:velocity-update-1}
    M_{n+1/2}^{(i)}=\left(e^{-s_n} + s_n\right) M^{\prime(i)}_{n} +  A_{n}^{(i)}\Delta t + \sqrt{D_n \cfrac{\tau_n}{2} ( 1 - e^{-2s_n}) } \xi^{(i)} , 
\end{eqnarray}
and the final velocity update
\begin{eqnarray}
     M_{n+1}^{(i)} &=& {\varepsilon_n} M_{n+1/2}^{(i)} + \hat{\mathbb{E}}^{(N)}[ M_{n} | X_n=x], \label{eq:velocity-update-2}
\end{eqnarray}
where $\varepsilon_n$ is the energy scaling factor
\begin{eqnarray}
    \varepsilon_n =  \cfrac{\hat{\mathbb{E}}^{(N)}\left[ \| M_{n+1/2} \|_2 ^2 \ \big| X_n = x\right]}{\hat{\mathbb{E}}^{(N)}\left[ \| M_{n} \|_2 ^2\ \big| \ X_n =  x \right]} \ .
\end{eqnarray}
Note that in the case of frozen coefficients, the presented scheme ensures conservation of momentum and energy during the velocity update and becomes effectively the Euler-Maruyama scheme as $s_n\to 0$. 
\noindent Particles positions are simply updated by
\begin{eqnarray} \label{eq:position-update}
    X_{n+1}^{(i)} = X_{n}^{(i)}+M_{n+1}^{(i)} \Delta t.
\end{eqnarray}

\subsection{Summary of algorithm} \label{subsection:final-algorithm}
\noindent We conclude the section with an overview on the full algorithm, whose steps are shown in Algorithm~\ref{alg:fefp}. The physical domain is discretized with a computational grid and at every iteration of the algorithm, particles are mapped to a computational cell according to their position. Since moments are evaluated conditionally on the position, they are computed cell-wise, as described in step \ref{algline:3}. 
Therefore, two main data structures are needed: one for the grid containing information on the moment estimates and drift coefficients, and another one for the particles, with every particle storing its position and velocity. Step \ref{algline:5} is used for post processing and typically requires a separate data structure.
\begin{algorithm}
\caption{Full FE-FP time loop}\label{alg:fefp}
\begin{algorithmic}[1]
\FOR {$i = 0, \Delta t, 2\Delta t, \dots, n \Delta t$}
    \FOR {$j=1,\dots,N_p$, $N_p=$ total number of particles,}
    \STATE Update particle position using \eqref{eq:position-update} \label{algline:1}
    \ENDFOR
    \STATE Apply boundary conditions
    \STATE Assign particles to computational cells
    \FOR{$j=1,\dots,N_c$, $N_c=$ total number of cells,}
    \STATE Update $c^{(4)} $ using \eqref{def:stab-coeff} \label{algline:2}
    \STATE Compute relevant moments using \eqref{def:montecarlo-integration}  \label{algline:3}
    \STATE Solve linear systems \eqref{eq:c-tilde} and \eqref{eq:linear-system}
    \STATE Compute $\{ \tilde{c}^{(1)} , c^{(2)} , c^{(3)}  \}$ through relation \eqref{eq:coeff-relation} \label{agline:3.1}
    \IF{$i\Delta t \geq $  steady state}
    \STATE Time average macroscopic quantities for output \label{algline:5}
    \ENDIF
    \ENDFOR
 \FOR {$j=1,\dots,N_p$,}
    \STATE Construct $A_n $ and $D_n $ based on particle location
    \STATE Update particle velocities according to \eqref{eq:velocity-update-1} and \eqref{eq:velocity-update-2} \label{algline:4} 
    \ENDFOR    
\ENDFOR 
\end{algorithmic}
\end{algorithm}

\section{Numerical results} \label{sec:results}
\begin{figure}
    \centering
    \begin{subfigure}[b]{0.49\textwidth}
        \centering
        \includegraphics[width=\textwidth]{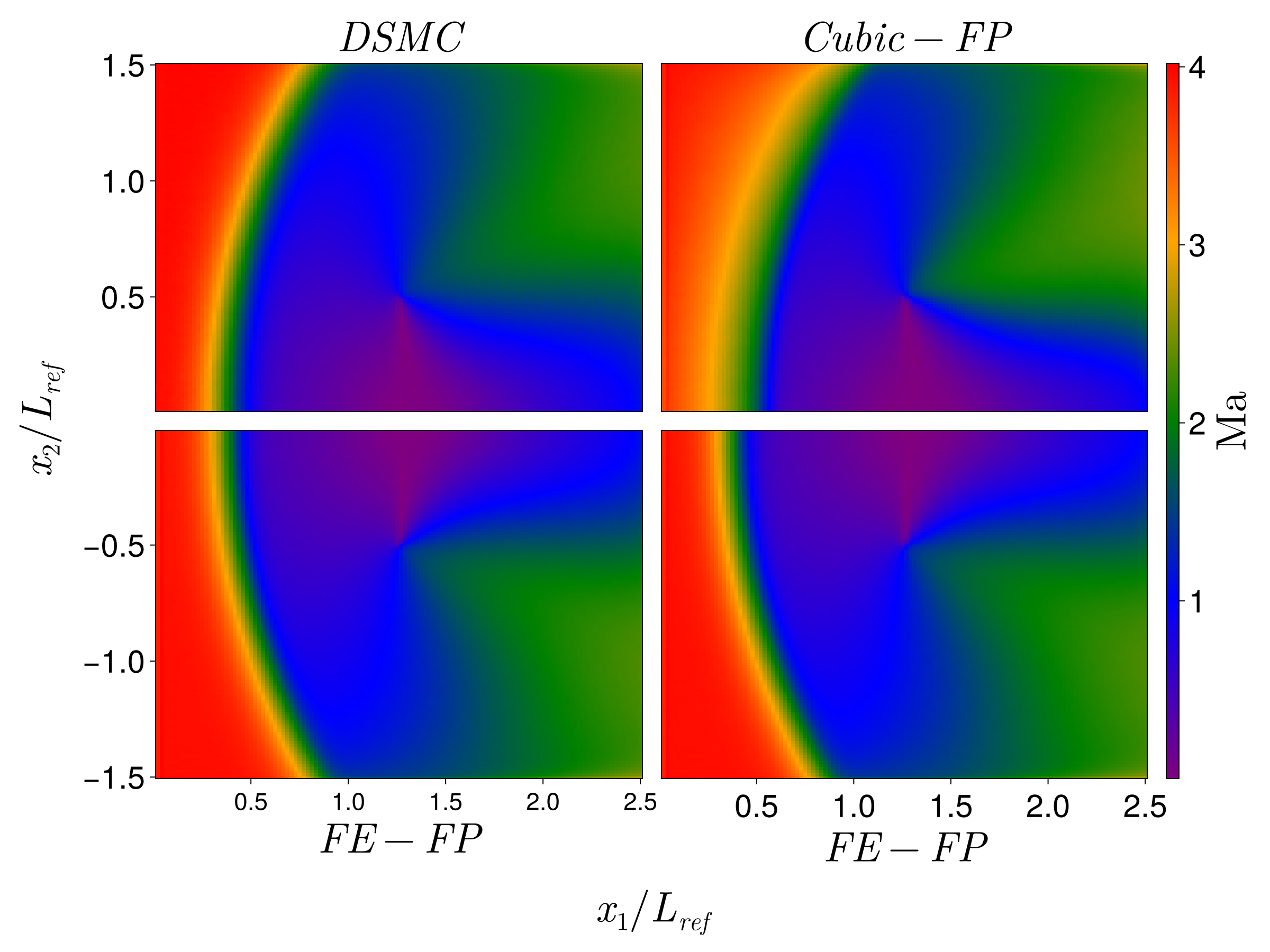}
        \caption{}
    \label{fig:mach-4}
    \end{subfigure}
    \hfill
    \begin{subfigure}[b]{0.49\textwidth}
        \centering
        \includegraphics[width=\textwidth]{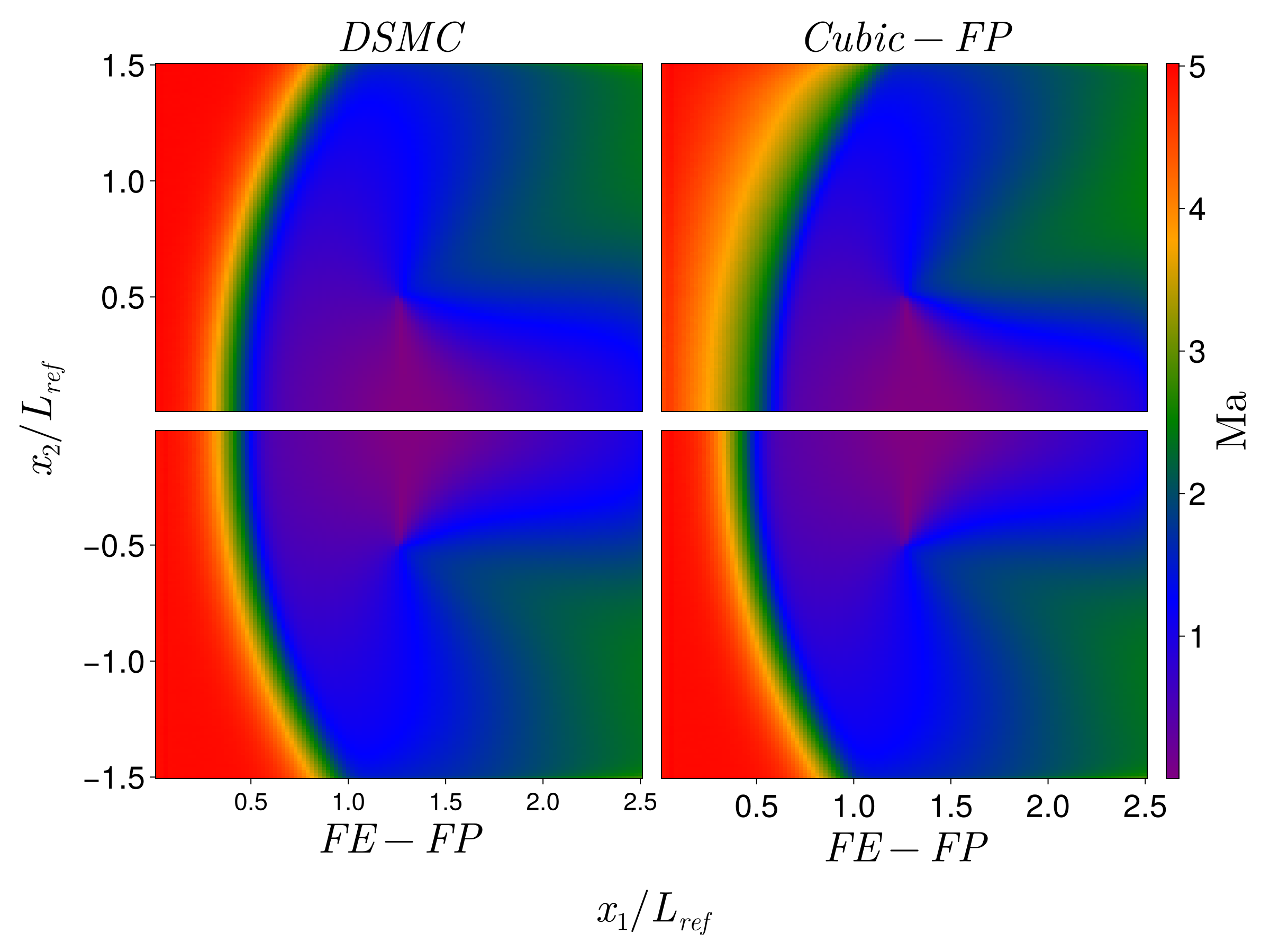}
        \caption{}
    \label{fig:mach-4}  
    \end{subfigure}
    \hfill
    \begin{subfigure}[b]{0.49\textwidth}
        \centering
        \includegraphics[width=\textwidth]{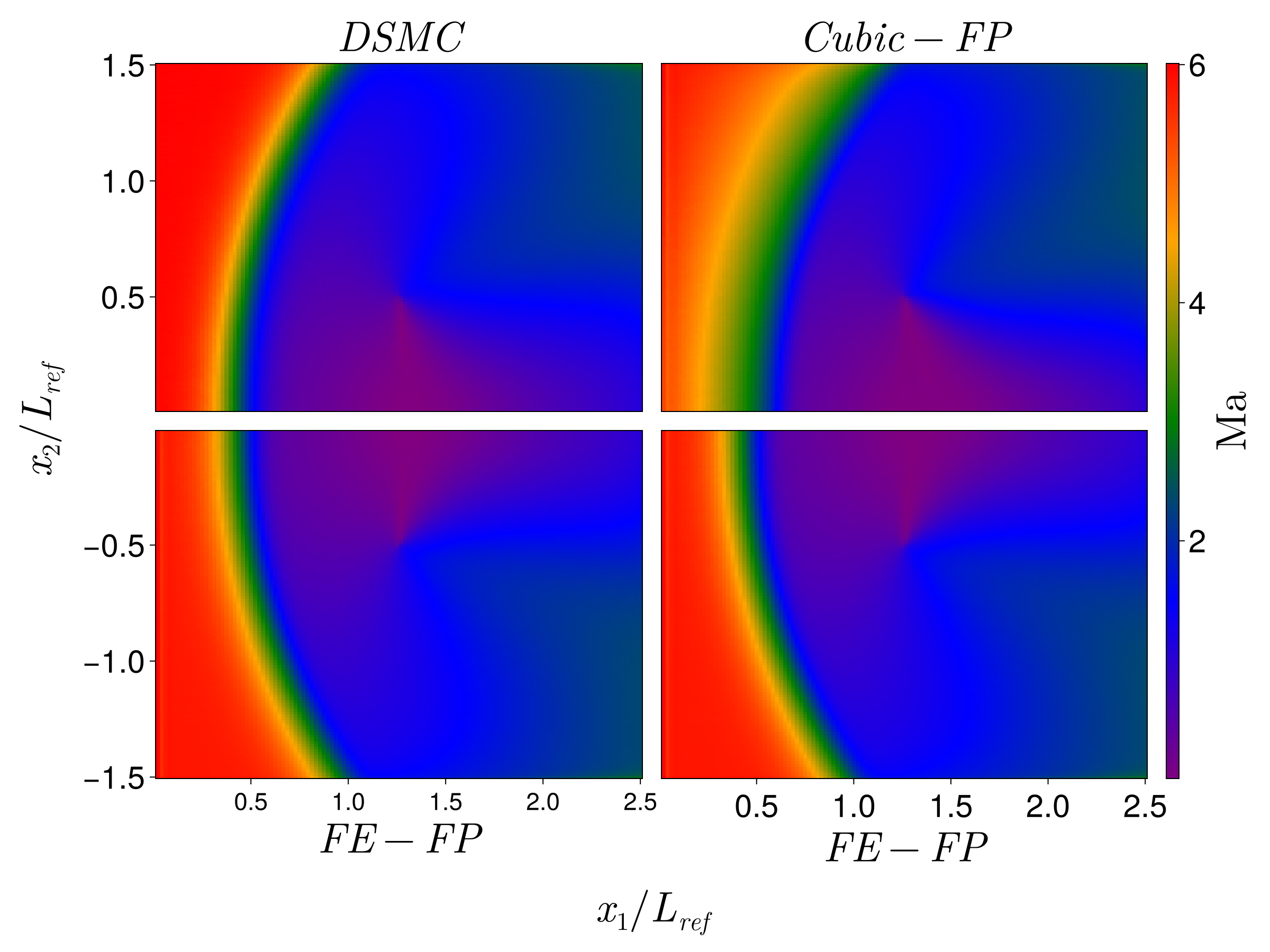}
        \caption{}
        \label{fig:mach-6}
    \end{subfigure}
    \hfill
    \begin{subfigure}[b]{0.49\textwidth}
        \centering
        \includegraphics[width=\textwidth]{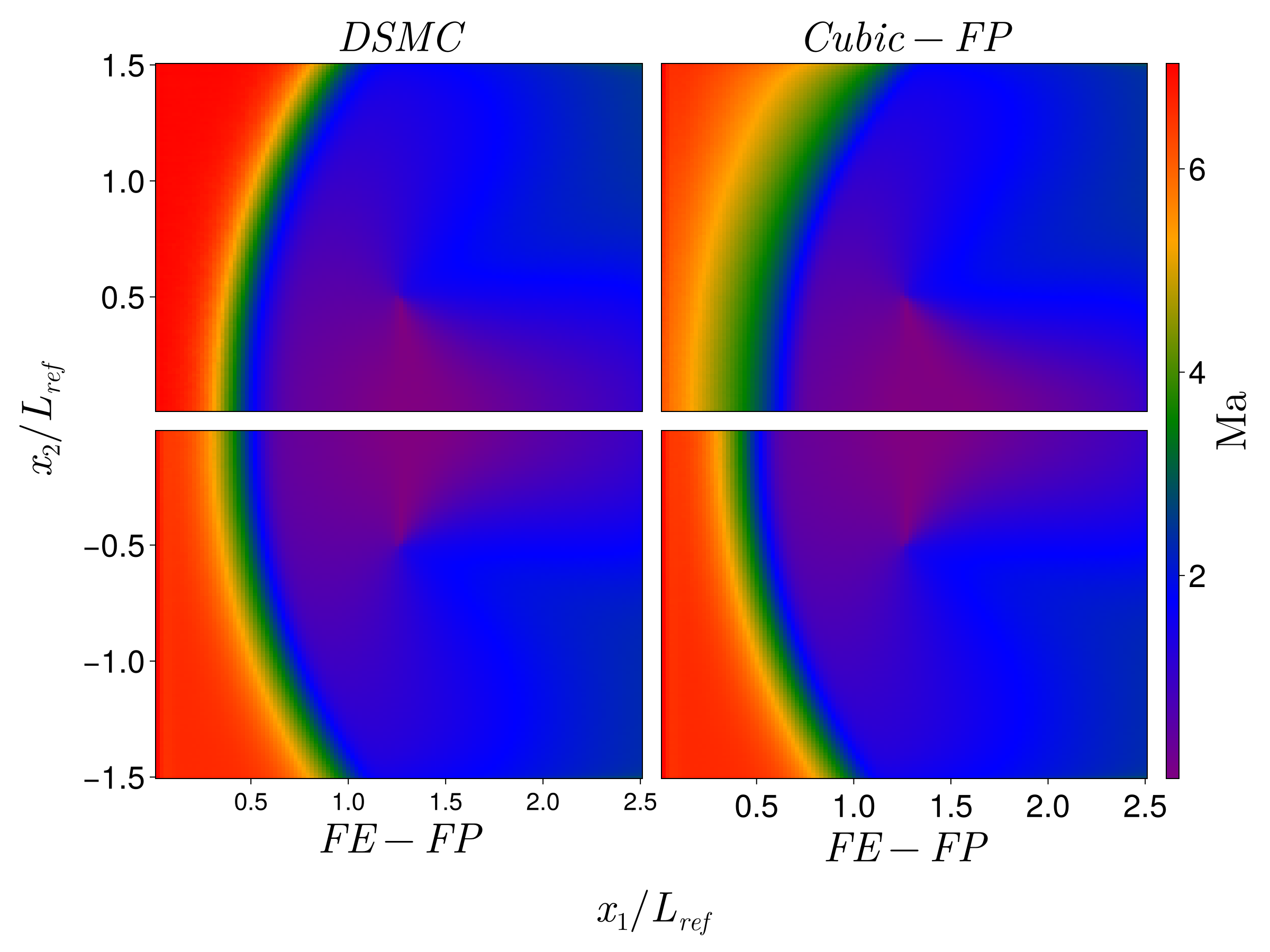}
        \caption{}
        \label{fig:mach-7}
    \end{subfigure}
    \caption{Mach contours for a flow of argon over a vertical plate; Results for (a) $\textrm{Ma}=4$, (b) $\textrm{Ma}=5$ , (c) $\textrm{Ma}=6$, and (d) $\textrm{Ma}=7$; DSMC and cubic-drift FP results are shown on top and compared with respect to the FE-FP results at the bottom.}
    \label{fig:temp-contour}
\end{figure}

\begin{figure}
    \centering
    \begin{subfigure}[h]{\textwidth}
        \begin{subfigure}[h]{0.46\textwidth}
            \centering
            \includegraphics[width=0.8\textwidth]{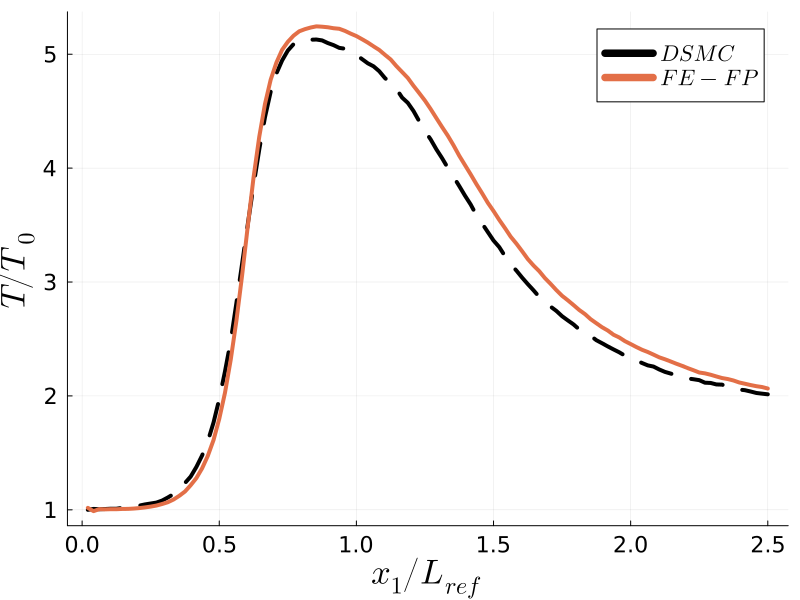}
            \caption{}
        \end{subfigure}
        \hfill
        \begin{subfigure}[h]{0.46\textwidth}
            \centering
            \includegraphics[width=0.8\textwidth]{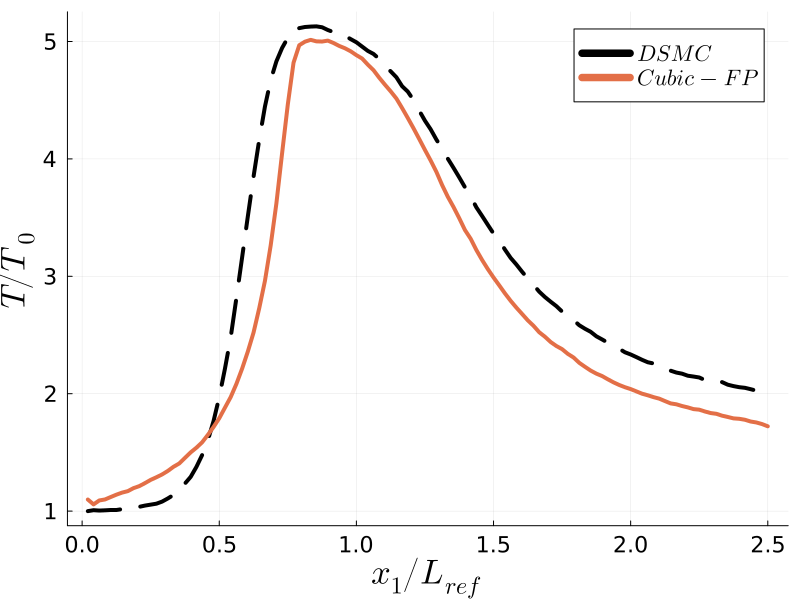}
            \caption{}
        \end{subfigure}
        \label{fig:temperature:4}
    \end{subfigure}
    \hfill
    \begin{subfigure}[h]{\textwidth}
        \begin{subfigure}[h]{0.46\textwidth}
            \centering
            \includegraphics[width=0.8\textwidth]{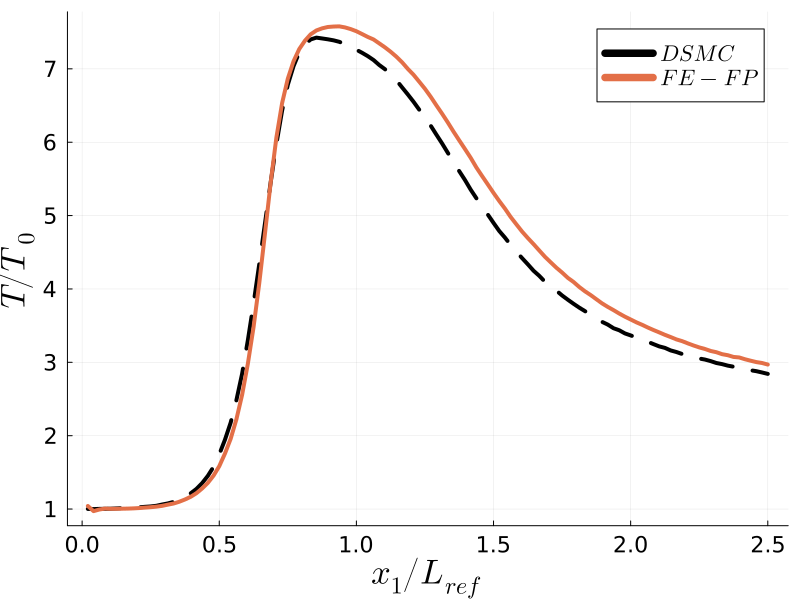}
            \caption{}
        \end{subfigure}
        \hfill
        \begin{subfigure}[h]{0.46\textwidth}
            \centering
            \includegraphics[width=0.8\textwidth]{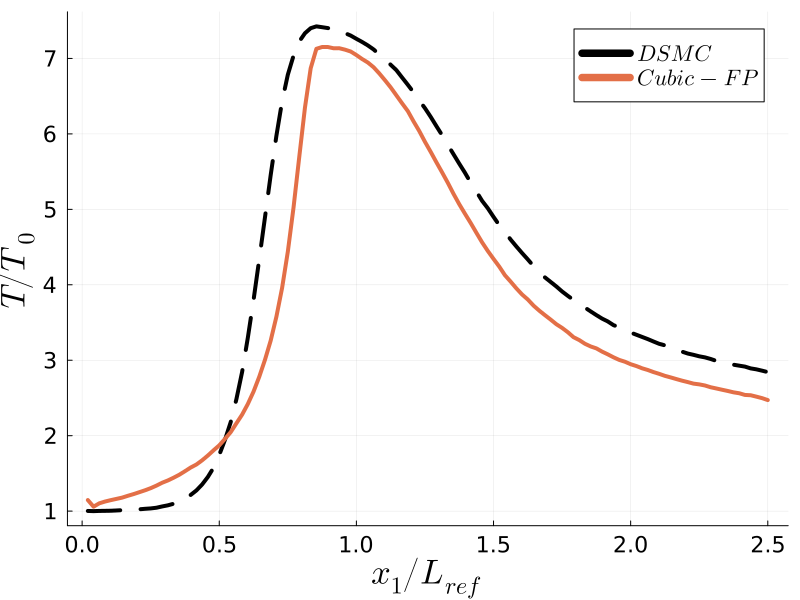}
            \caption{}
        \end{subfigure}
        \label{fig:temperature:5}
    \end{subfigure}
    \hfill
    \begin{subfigure}[h]{\textwidth}
        \begin{subfigure}[h]{0.46\textwidth}
            \centering
            \includegraphics[width=0.8\textwidth]{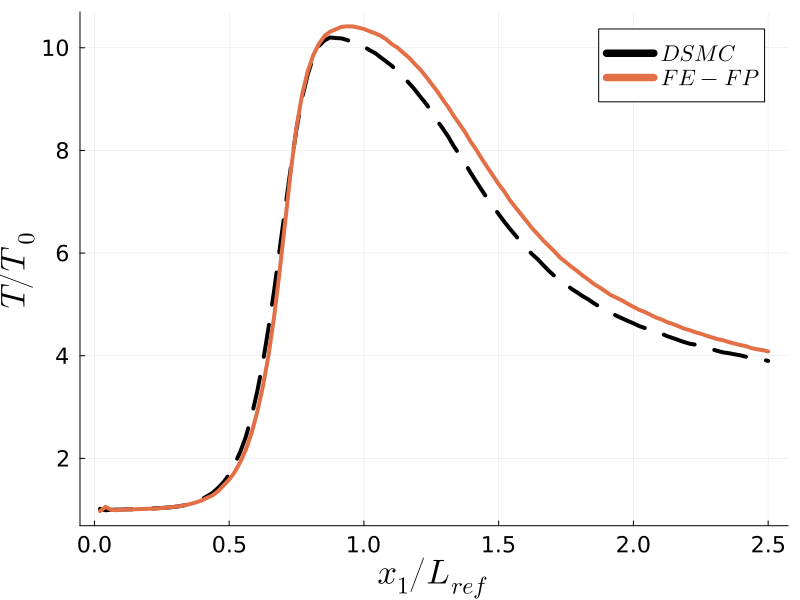}
            \caption{}
        \end{subfigure}
        \hfill
        \begin{subfigure}[h]{0.46\textwidth}
            \centering
            \includegraphics[width=0.8\textwidth]{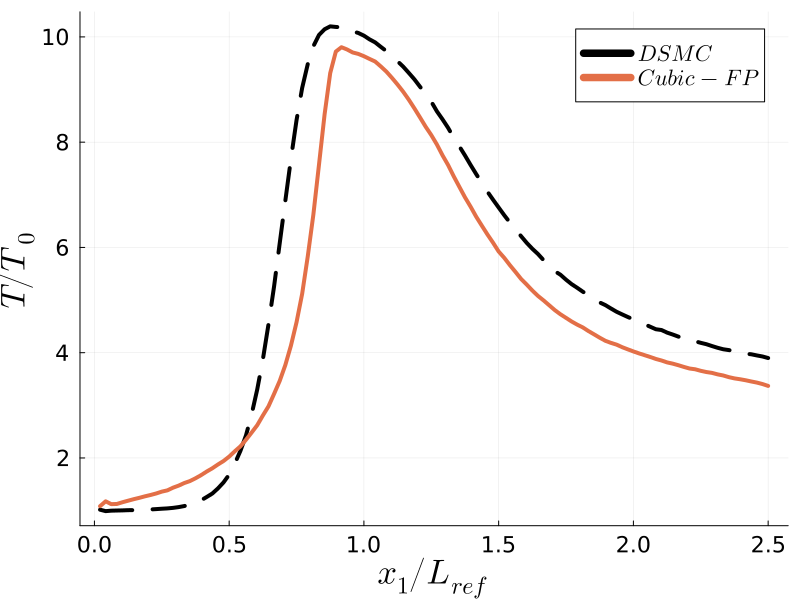}
            \caption{}
        \end{subfigure}
        \label{fig:temperature:6}
    \end{subfigure}
    \hfill
    \begin{subfigure}[h]{\textwidth}
        \begin{subfigure}[h]{0.46\textwidth}
            \centering
            \includegraphics[width=.8\textwidth]{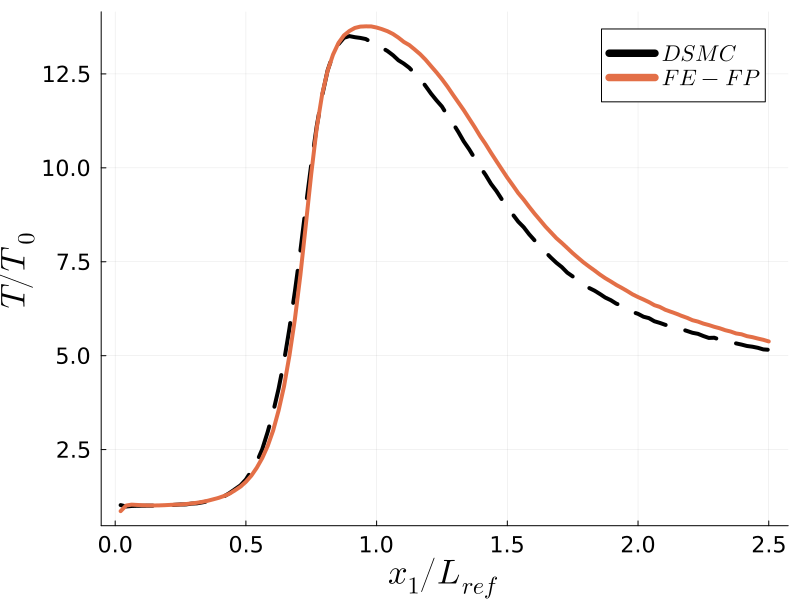}
            \caption{}
        \end{subfigure}
        \hfill
        \begin{subfigure}[h]{0.46\textwidth}
            \centering
            \includegraphics[width=.8\textwidth]{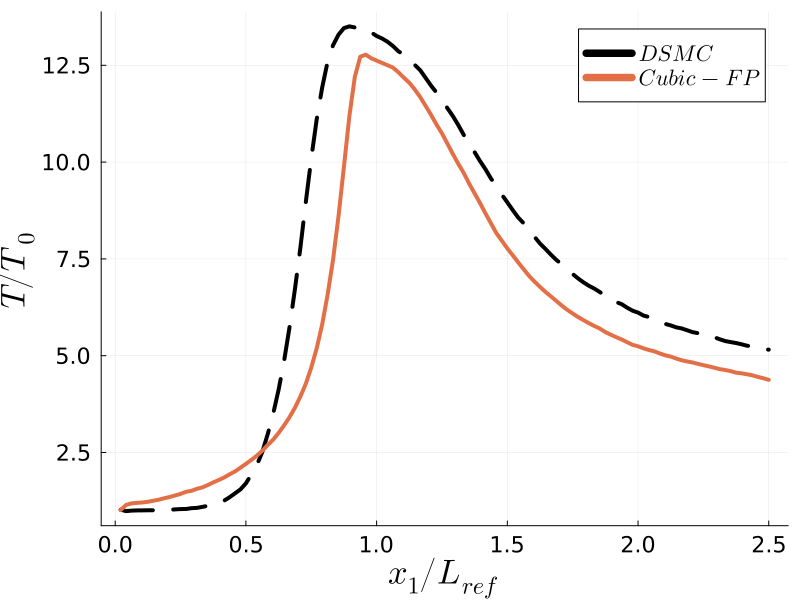}
            \caption{}
        \end{subfigure}
        \label{fig:temperature:7}
    \end{subfigure}
    \caption{Normalized temperature profiles along $x_2=1.875L_{ref}$; (a)-(b) $\textrm{Ma}=4$, (c)-(d) $\textrm{Ma}=5$, (e)-(f) $\textrm{Ma}=6$, and (g)-(h) $\textrm{Ma}=7$; The left column compares the DSMC reference solution with the FE-FP model, while the right column shows comparisons with the cubic-drift FP model.}
    \label{fig:temp-profile}
\end{figure}

\noindent We present numerical results of FE-FP in recovering the shock profiles, where comparison with respect to the cubic-drift FP model \cite{gorjiFokkerPlanckDSMC2015}, and benchmark DSMC \cite{birdOneDimensionalSteadyFlows1994} is provided. The setting chosen for validation of the FE-FP model is a supersonic flow of argon (modeled as hard sphere gas) over an infinitely thin plate \cite{gorjiEntropicFokkerPlanckKinetic2021}. \\ \ \\ \noindent
A fully diffusive isothermal plate of reference length $L_{ref}$ and temperature $T_{ref}$ lies at the center of a rectangular domain of size $L_{x_1}\times L_{x_2}$ with $L_{x_1}=2.5L_{ref}$ and $L_{x_2} = 3 L_{ref}$, where the gas with initial number density $n_0$ is at equilibrium with temperature $T_0=T_{ref}$, initially. Calculating the reference time scale $\tau_{ref} = 2\mu_0/p_0$ from the initial viscosity $\mu_0=\mu(T_0)$ and pressure $p_0=n_0 k_BT_0$. We choose the Knudsen number $Kn=0.14$ for the length $L_{ref}/2$ and the hard-sphere mean-free-path $\lambda$. The supersonic flow enters through the left boundary resulting in a shock around the plate. For hard sphere gas, the viscosity law is $\mu(T) = \mu(T_0) \sqrt{T/T_0}$. \\ \ \\
The computational grid is comprised of $120\times120$  cells for discretization of physical space with the horizontal and vertical coordinates of $x_1$ and $x_2$, respectively. The initial number of computational particles is $3\times10^6$. The timestep is chosen fine enough such that the benchmark DSMC calculations are accurate, leading to $\Delta t= 1.3876\times 10^{-2}\tau_{ref}$. To reach steady state, the simulation is run for $2,000$ time-steps, and additional $5,000$ time-steps are used to average the macroscopic quantities. 
Because of the $x_2$-axis symmetry of the problem, only the upper half of the domain is simulated where symmetric boundary condition is applied on $x_2=0$ line. For the open-boundaries including both outflow and inflow boundary conditions, we use the technique described by Bird in \cite{birdOneDimensionalSteadyFlows1994}, while particles leaving the domain are deleted from the simulation. The simulation is repeated for upstream Mach number $\textrm{Ma}=[4,\ 5, \ 6, \ 7]$. The FE-FP model setting is the one described in Sec.~\ref{sec:num}. Moments matching up until heat-fluxes are imposed for hard sphere production terms as described in \cite{guptaAutomatedBoltzmannCollision2012}. 
\\ \ \\
\noindent
The shock profile at steady state, visualized through Mach contours, is shown in Fig.~\ref{fig:temp-contour}. Very good agreement on the shock profile is shown between the DSMC reference solution and FE-FP for the whole range of Mach numbers. On the other hand, the cubic-drift FP shows a diffusive behavior that progressively gets more evident as the Mach number increases. For better visualization, in Fig.~\ref{fig:temp-profile} we depict the temperature profiles with respect to $x_1$ by slicing the domain at $x_2=1.875 L_{ref}$, slightly above the plate. The agreement with the DSMC solution is further confirmed there. \\ \ \\
Despite the overall close agreement between FE–FP and DSMC results, we note a visible deviation in the shock downstream region. These discrepancies are expected in zones characterized by large local Knudsen numbers, where non-equilibrium effects can extend beyond those captured by moment matching up to the heat-fluxes. We anticipate that augmenting the FE–FP framework with higher-order moment constraints would help diminish this gap and further improve the accuracy of shock-structure predictions in such regimes.

\section{Conclusion and final remarks}
\noindent Kinetic models have attained growing attention as an active area of research, bridging interacting particle systems with macroscopic continuum descriptions in fields ranging from rarefied gas dynamics to plasma physics and data-driven modeling. Among these, the FP equation offers a particularly attractive formulation, arising as the diffusive approximation of the Boltzmann equation in regimes where collisions are overwhelming. Yet, even in the classical setting of a monatomic dilute gas, the rigorous derivation of such approximations remains an open challenge.
\\ \ \\
Our study addresses this problem by introducing the FE-FP equation. The model brings together a consistent moment projection of the Boltzmann equation with an entropy-decay mechanism with favorable structure. The formulation is both analytically tractable and computationally affordable, requiring only marginal computational overhead compared to standard FP models. We provide theoretical justifications for the basic aspects of the model, including the entropy decay and the well-posedness of linear system underlying the model coefficients. The framework is also validated by numerical results on shock waves  confirming the role of correct entropy production in shock profiles, showing close agreement with respect to the high-fidelity solution of the Boltzmann equation. 
\\ \ \\
While the present study emphasizes on the model development, relying on standard numerical solvers for the resulting SDE systems, the FE–FP formulation is designed as a flexible framework upon which efficient, high-order integration schemes can be built. Such developments, aimed at improving convergence rate of resulting particle schemes, will be pursued in future work.
\\ \ \\
The FE-FP model supports arbitrary degrees of freedom in moment matching, makes it amenable to \textit{a posteriori} error estimates. Following this path, the FE-FP model can be employed to build cascades of FP models, extending the application range of FP models to simulations of multiscale problems. More broadly, the proposed framework provides a path to building SDE models constrained by both moment conditions and entropy laws, with potential applications extending beyond kinetic theory to areas such as diffusion-based generative modeling and the Schrödinger bridge problem \cite{gottwald2025stable}.

\section*{Acknowledgments}
\noindent V.M. gratefully acknowledges Dr. Ivan Lunati for providing laboratory infrastructure, computational resources, and insightful feedback, and Prof. Fabio Nobile for his valuable comments on the study.

\section*{Reproducibility of computational results}
 \noindent Code and data that allow readers to reproduce the results in this paper are available at \url{https://codeberg.org/vmontanaro/StochasticFP}.

\appendix

\section{Positive defineteness of $R[f]$} \label{appendix:a}
\noindent Consider an arbitrary non-zero $z \in \mathbb{R}^l$. We have
\begin{eqnarray} \label{eq:matrix-pd}
    \sum_{\alpha,\beta \in \hat{I}_{2k}^3 } z_\alpha R_{\alpha\beta}[f]z_\beta &=&\left\langle \left( \sum_{\alpha \in \hat{I}_{2k}^3} z_\alpha \nabla_v H_{\alpha}(v) \right) \cdot \left( \sum_{\beta \in \hat{I}_{2k}^3} z_\beta \nabla_v H_{\beta}(v) \right), f\right\rangle \nonumber \\
    &=& \left\langle \left( \sum_{\alpha \in \hat{I}_{2k}^3} z_\alpha \nabla_v H_{\alpha}(v) \right)^2, f \right \rangle \ .
\end{eqnarray}
\sloppy The right-hand side is clearly non-negative, so we only need to check whether $ \sum_{\alpha \in \hat{I}_{2k}^3} z_\alpha \nabla_v H_{\alpha}(v)  = 0$. This is prevented due to the linear independence of  $\{H_\alpha(v)\}_{\alpha \in \hat{I}^{3}_{2k}}$ polynomials (by construct). Therefore the summation is strictly positive for all $z \in \mathbb{R}^l$. Besides the fact that $R[f]$ is symmetric, this concludes the proof that $R[f]$ is symmetric positive definite.
\section{Three-dimensional linear system matrix} \label{appendix:c}
\noindent Here we provide the elements of matrices and vectors for the linear system corresponding to the model coefficients of Sec.~\ref{subsec:model-3d}. To simplify notation, we drop the dependency on $f$ and leave it implied for the rest of the appendix.  Moreover, we use the naming $R,Q$ and $b$ as defined in Sec.~\ref{subsec:linear-system}.
We have 
\begin{eqnarray}
    R = \begin{pmatrix}
        B_{6 \times6} & C_{6 \times 3} & E_{6\times1} \\
        F_{3\times6} & J_{3\times3} & Z_{3\times1} \\
        E_{1\times6} & Z_{1\times3} & W_{1\times1} \\
    \end{pmatrix} \ ,
\end{eqnarray}
where
\begin{eqnarray}
    &&B = \begin{pmatrix}
        2\tilde{\mathrm{m}}^{(0)}_{11} & 2\tilde{\mathrm{m}}^{(0)}_{12} & 2\tilde{\mathrm{m}}^{(0)}_{13} & 0 & 0 & 0 \\
        \tilde{\mathrm{m}}^{(0)}_{12} & \tilde{\mathrm{m}}^{(0)}_{11} + \tilde{\mathrm{m}}^{(0)}_{22} & \tilde{\mathrm{m}}^{(0)}_{23} & \tilde{\mathrm{m}}^{(0)}_{12} & \tilde{\mathrm{m}}^{(0)}_{13} & 0 \\
        \tilde{\mathrm{m}}^{(0)}_{13} & \tilde{\mathrm{m}}^{(0)}_{23} & \tilde{\mathrm{m}}^{(0)}_{11} + \tilde{\mathrm{m}}^{(0)}_{33} & 0 & \tilde{\mathrm{m}}^{(0)}_{12} & \tilde{\mathrm{m}}^{(0)}_{13} \\
        0 & 2\tilde{\mathrm{m}}^{(0)}_{12} & 0 & 2\tilde{\mathrm{m}}^{(0)}_{22} & 2\tilde{\mathrm{m}}^{(0)}_{23} & 0 \\
        0 & \tilde{\mathrm{m}}^{(0)}_{13} & \tilde{\mathrm{m}}^{(0)}_{12} & \tilde{\mathrm{m}}^{(0)}_{23} & \tilde{\mathrm{m}}^{(0)}_{22} + \tilde{\mathrm{m}}^{(0)}_{33} & \tilde{\mathrm{m}}^{(0)}_{23} \\
        0 & 0 & 2\tilde{\mathrm{m}}^{(0)}_{13} & 0 & \tilde{\mathrm{m}}^{(0)}_{23} & 2\tilde{\mathrm{m}}^{(0)}_{33}
        
    \end{pmatrix} \ ,
\end{eqnarray}
\begin{eqnarray}
    C = \begin{pmatrix}
        2\tilde{\mathrm{m}}^{(2)}_{1} +4\tilde{\mathrm{m}}^{(0)}_{111} & 4\tilde{\mathrm{m}}^{(0)}_{112} & 4\tilde{\mathrm{m}}^{(0)}_{113} \\
        \tilde{\mathrm{m}}^{(2)}_{1} + 4\tilde{\mathrm{m}}^{(0)}_{112} & \tilde{\mathrm{m}}^{(2)}_{2} 4\tilde{\mathrm{m}}^{(0)}_{122} & 4\tilde{\mathrm{m}}^{(0)}_{123} \\
        \tilde{\mathrm{m}}^{(2)}_{1} + 4\tilde{\mathrm{m}}^{(0)}_{113} & 4\tilde{\mathrm{m}}^{(0)}_{123} & \tilde{\mathrm{m}}^{(2)}_{3} + 4\tilde{\mathrm{m}}^{(0)}_{133} \\
        4\tilde{\mathrm{m}}^{(0)}_{122} & 2\tilde{\mathrm{m}}^{(2)}_{2} + 4\tilde{\mathrm{m}}^{(0)}_{222} & 4\tilde{\mathrm{m}}^{(0)}_{223} \\
        4\tilde{\mathrm{m}}^{(0)}_{123} & \tilde{\mathrm{m}}^{(2)}_{2} + 4\tilde{\mathrm{m}}^{(0)}_{223} & \tilde{\mathrm{m}}^{(2)}_{3} + 4\tilde{\mathrm{m}}^{(0)}_{233} \\
        4\tilde{\mathrm{m}}^{(0)}_{133} & 4\tilde{\mathrm{m}}^{(0)}_{233} & 2\tilde{\mathrm{m}}^{(2)}_{3} + 4\tilde{\mathrm{m}}^{(0)}_{333}
    \end{pmatrix} \ ,
\end{eqnarray}
\begin{eqnarray}
    E = 6\begin{pmatrix}
        \tilde{\mathrm{m}}^{(2)}_{11} & \tilde{\mathrm{m}}^{(2)}_{12} & \tilde{\mathrm{m}}^{(2)}_{13} & \tilde{\mathrm{m}}^{(2)}_{22} & \tilde{\mathrm{m}}^{(2)}_{23} & \tilde{\mathrm{m}}^{(2)}_{33} 
    \end{pmatrix} \ ,
\end{eqnarray}
\begin{eqnarray}
    F &=& \begin{pmatrix}
       2\tilde{\mathrm{m}}^{(0)}_{111} & 4\tilde{\mathrm{m}}^{(0)}_{112} & 4\tilde{\mathrm{m}}^{(0)}_{113} & 2\tilde{\mathrm{m}}^{(0)}_{122} & 4\tilde{\mathrm{m}}^{(0)}_{123} & 2\tilde{\mathrm{m}}^{(0)}_{133} \\ 
        2\tilde{\mathrm{m}}^{(0)}_{112} &  4\tilde{\mathrm{m}}^{(0)}_{122} & 4\tilde{\mathrm{m}}^{(0)}_{123} &  2\tilde{\mathrm{m}}^{(0)}_{222} &  4\tilde{\mathrm{m}}^{(0)}_{223} & 2\tilde{\mathrm{m}}^{(0)}_{233} \\
        2\tilde{\mathrm{m}}^{(0)}_{113} & 4\tilde{\mathrm{m}}^{(0)}_{123} & 4\tilde{\mathrm{m}}^{(0)}_{133} & 2\tilde{\mathrm{m}}^{(0)}_{223} & 4\tilde{\mathrm{m}}^{(0)}_{233} &  2\tilde{\mathrm{m}}^{(0)}_{333}
    \end{pmatrix} \nonumber \\
    &+& \begin{pmatrix}
         \tilde{\mathrm{m}}^{(2)}_1 &  \tilde{\mathrm{m}}^{(2)}_2 &  \tilde{\mathrm{m}}^{(2)}_3 & 0 & 0 & 0 \\
         0 & \tilde{\mathrm{m}}^{(2)}_1  & 0 & \tilde{\mathrm{m}}^{(2)}_2 & \tilde{\mathrm{m}}^{(2)}_3 & 0 \\
         0 & 0 & \tilde{\mathrm{m}}^{(2)}_1 & 0 & \tilde{\mathrm{m}}^{(2)}_2 & \tilde{\mathrm{m}}^{(2)}_3
    \end{pmatrix} \ ,
\end{eqnarray}
\sloppy
\begin{eqnarray}
    J &=& \begin{pmatrix}
        8\tilde{\mathrm{m}}^{(2)}_{11} -4\tilde{\mathrm{m}}^{(0)}_{11}\tilde{\mathrm{m}}^{(2)} & 8\tilde{\mathrm{m}}^{(2)}_{12} -4\tilde{\mathrm{m}}^{(0)}_{12}\tilde{\mathrm{m}}^{(2)} & 8\tilde{\mathrm{m}}^{(2)}_{13} -4\tilde{\mathrm{m}}^{(0)}_{13}\tilde{\mathrm{m}}^{(2)} \\
        
        8\tilde{\mathrm{m}}^{(2)}_{12} -4\tilde{\mathrm{m}}^{(0)}_{12}\tilde{\mathrm{m}}^{(2)} & 8\tilde{\mathrm{m}}^{(2)}_{22} -4\tilde{\mathrm{m}}^{(0)}_{22}\tilde{\mathrm{m}}^{(2)} & 8\tilde{\mathrm{m}}^{(2)}_{23} -4\tilde{\mathrm{m}}^{(0)}_{23}\tilde{\mathrm{m}}^{(2)} \\
        
        8\tilde{\mathrm{m}}^{(2)}_{13} -4\tilde{\mathrm{m}}^{(0)}_{13}\tilde{\mathrm{m}}^{(2)} & 8\tilde{\mathrm{m}}^{(2)}_{23} -4\tilde{\mathrm{m}}^{(0)}_{23}\tilde{\mathrm{m}}^{(2)} & 8\tilde{\mathrm{m}}^{(2)}_{33} -4\tilde{\mathrm{m}}^{(0)}_{33}\tilde{\mathrm{m}}^{(2)} \\
    \end{pmatrix}  \nonumber \\
    &-& 4\begin{pmatrix}
        (\tilde{\mathrm{m}}^{(0)}_{11})^2 + (\tilde{\mathrm{m}}^{(0)}_{12})^2 + (\tilde{\mathrm{m}}^{(0)}_{13})^2  &
        0 & 0 \\
        0 & (\tilde{\mathrm{m}}^{(0)}_{12})^2 + (\tilde{\mathrm{m}}^{(0)}_{22})^2 + (\tilde{\mathrm{m}}^{(0)}_{23})^2  & 0 \\
        0 & 0 & (\tilde{\mathrm{m}}^{(0)}_{13})^2 + (\tilde{\mathrm{m}}^{(0)}_{23})^2 + (\tilde{\mathrm{m}}^{(0)}_{33})^2 
    \end{pmatrix}   \nonumber \\
    &-& 4\begin{pmatrix}
        0 & \tilde{\mathrm{m}}^{(0)}_{11}\tilde{\mathrm{m}}^{(0)}_{12} & \tilde{\mathrm{m}}^{(0)}_{11}\tilde{\mathrm{m}}^{(0)}_{13} \\
        \tilde{\mathrm{m}}^{(0)}_{11}\tilde{\mathrm{m}}^{(0)}_{12} & 0 & \tilde{\mathrm{m}}^{(0)}_{12}\tilde{\mathrm{m}}^{(0)}_{13} \\
        \tilde{\mathrm{m}}^{(0)}_{11}\tilde{\mathrm{m}}^{(0)}_{13} & \tilde{\mathrm{m}}^{(0)}_{11}\tilde{\mathrm{m}}^{(0)}_{13} & 0
    \end{pmatrix} - 4 \begin{pmatrix}
     0 & \tilde{\mathrm{m}}^{(0)}_{12}\tilde{\mathrm{m}}^{(0)}_{22} & \tilde{\mathrm{m}}^{(0)}_{12}\tilde{\mathrm{m}}^{(0)}_{23} \\
        \tilde{\mathrm{m}}^{(0)}_{12}\tilde{\mathrm{m}}^{(0)}_{22} & 0 & \tilde{\mathrm{m}}^{(0)}_{22}\tilde{\mathrm{m}}^{(0)}_{23} \\
        \tilde{\mathrm{m}}^{(0)}_{12}\tilde{\mathrm{m}}^{(0)}_{23} & \tilde{\mathrm{m}}^{(0)}_{12}\tilde{\mathrm{m}}^{(0)}_{23} & 0
    \end{pmatrix}\nonumber  \\
    &-& 4 \begin{pmatrix}
             0 & \tilde{\mathrm{m}}^{(0)}_{13}\tilde{\mathrm{m}}^{(0)}_{23} & \tilde{\mathrm{m}}^{(0)}_{13}\tilde{\mathrm{m}}^{(0)}_{33} \\
        \tilde{\mathrm{m}}^{(0)}_{13}\tilde{\mathrm{m}}^{(0)}_{23} & 0 & \tilde{\mathrm{m}}^{(0)}_{23}\tilde{\mathrm{m}}^{(0)}_{33} \\
        \tilde{\mathrm{m}}^{(0)}_{13}\tilde{\mathrm{m}}^{(0)}_{33} & \tilde{\mathrm{m}}^{(0)}_{13}\tilde{\mathrm{m}}^{(0)}_{33} & 0
    \end{pmatrix} + \begin{pmatrix}
        \tilde{\mathrm{m}}^{(4)} - (\tilde{\mathrm{m}}^{(2)})^2 & 0 & 0 \\
        0 & \tilde{\mathrm{m}}^{(4)} - (\tilde{\mathrm{m}}^{(2)})^2 & 0 \\
        0 & 0 & \tilde{\mathrm{m}}^{(4)} - (\tilde{\mathrm{m}}^{(2)})^2
    \end{pmatrix} \ , \nonumber \\
\end{eqnarray}
\begin{eqnarray}
    Z = 3 \begin{pmatrix}
        3\tilde{\mathrm{m}}^{(4)}_1 - 2(\tilde{\mathrm{m}}^{(2)}_1 \tilde{\mathrm{m}}^{(0)}_{11} + \tilde{\mathrm{m}}^{(2)}_2\tilde{\mathrm{m}}^{(0)}_{12} + \tilde{\mathrm{m}}^{(2)}_3\tilde{\mathrm{m}}^{(0)}_{13}) - \tilde{\mathrm{m}}^{(2)}_1\tilde{\mathrm{m}}^{(2)} \\
        3\tilde{\mathrm{m}}^{(4)}_2 - 2(\tilde{\mathrm{m}}^{(2)}_1 \tilde{\mathrm{m}}^{(0)}_{12} + \tilde{\mathrm{m}}^{(2)}_2\tilde{\mathrm{m}}^{(0)}_{22} + \tilde{\mathrm{m}}^{(2)}_3\tilde{\mathrm{m}}^{(0)}_{23}) - \tilde{\mathrm{m}}^{(2)}_1\tilde{\mathrm{m}}^{(2)} \\
        3\tilde{\mathrm{m}}^{(4)}_1 - 2(\tilde{\mathrm{m}}^{(2)}_1 \tilde{\mathrm{m}}^{(0)}_{13} + \tilde{\mathrm{m}}^{(2)}_2\tilde{\mathrm{m}}^{(0)}_{23} + \tilde{\mathrm{m}}^{(2)}_3\tilde{\mathrm{m}}^{(0)}_{33}) - \tilde{\mathrm{m}}^{(2)}_1\tilde{\mathrm{m}}^{(2)}
    \end{pmatrix} 
\end{eqnarray}
and
\begin{eqnarray}
    W = 6\ \left(\tilde{\mathrm{m}}^{(4)}_{11}+\tilde{\mathrm{m}}^{(4)}_{22}+\tilde{\mathrm{m}}^{(4)}_{33} - (\tilde{\mathrm{m}}^{(2)}_1)^2 - (\tilde{\mathrm{m}}^{(2)}_2)^2 - (\tilde{\mathrm{m}}^{(2)}_3)^2 \right) \ .
\end{eqnarray}
The vector corresponding to the entropy constraint is given by
\begin{eqnarray}
    Q = \begin{pmatrix}
        X_{1\times6} & Y_{1 \times3} & 5\tilde{\mathrm{m}}^{(2)} 
    \end{pmatrix} \ ,
\end{eqnarray}
where
\begin{eqnarray}
    X = \begin{pmatrix}
        1 & 0 & 0 & 1 & 0 & 1
    \end{pmatrix} \rho 
\end{eqnarray}
and
\begin{eqnarray}
    Y = 8 \begin{pmatrix}
        \tilde{\mathrm{m}}^{(0)}_1 & \tilde{\mathrm{m}}^{(0)}_2 & \tilde{\mathrm{m}}^{(0)}_3
    \end{pmatrix}\ .
\end{eqnarray}
The right-hand side can be cast into
\begin{eqnarray}
    b = \begin{pmatrix}
        P^{(0)}_{11} - D + c^{(4)}\tilde{\mathrm{m}}^{(4)}_{11}\\
        P^{(0)}_{12} - 2c^{(4)} \tilde{\mathrm{m}}^{(4)}_{12} \\
        P^{(0)}_{13} - 2c^{(4)} 2\tilde{\mathrm{m}}^{(4)}_{13} \\
        P^{(0)}_{12} - D - c^{(4)} \tilde{\mathrm{m}}^{(4)}_{22} \\
        P^{(0)}_{23} - 2c^{(4)} \tilde{\mathrm{m}}^{(4)}_{23} \\
        P^{(0)}_{33} - D - c^{(4)} \tilde{\mathrm{m}}^{(4)}_{33} \\
        P^{(2)}_{1} - c^{(4)} \tilde{\mathrm{m}}^{(6)}_{1} \\
        P^{(2)}_{2}- c^{(4)} \tilde{\mathrm{m}}^{(6)}_{2} \\
        P^{(3)}_{1} - c^{(4)}\tilde{\mathrm{m}}^{(6)}_{3} \\
        - c^{(4)} \tilde{\mathrm{m}}^{(4)}
    \end{pmatrix} \ ,
\end{eqnarray}
where $\{P^{(p)}_\alpha\}_{\alpha \in \hat{I}^3_{4}}$ are the production terms.

\newpage

\bibliographystyle{siamplain}
\bibliography{references}

\begin{thebibliography}{10}

\bibitem{alderPhaseTransitionHard1957}
{\sc B.~Alder and T.~Wainwright}, {\em Phase transition for a hard sphere system}, The Journal of Chemical Physics, 27 (1957), pp.~1208--1209, \url{https://doi.org/10.1063/1.1743957}.

\bibitem{585795}
{\sc E.~Arkilic, M.~Schmidt, and K.~Breuer}, {\em Gaseous slip flow in long microchannels}, Journal of Microelectromechanical Systems, 6 (1997), pp.~167--178, \url{https://doi.org/10.1109/84.585795}.

\bibitem{arnoldLargetimeBehaviorNonsymmetric2008}
{\sc A.~Arnold, E.~Carlen, and Q.~Ju}, {\em Large-time behavior of non-symmetric {{Fokker-Planck}} type equations}, Communications on Stochastic Analysis, 2 (2008), \url{https://doi.org/10.31390/cosa.2.1.11}.

\bibitem{arnoldCONVEXSOBOLEVINEQUALITIES2001}
{\sc A.~Arnold, P.~Markowich, G.~Toscani, and A.~Unterreiter}, {\em {{On Convex Sobolev Inequalities and the Rate of Convergence to Equilibrium for Fokker-Planck Type Equations}}}, Communications in Partial Differential Equations, 26 (2001), pp.~43--100, \url{https://doi.org/10.1081/PDE-100002246}.

\bibitem{basovModelingPolyatomicGases2024}
{\sc L.~Basov and M.~Grabe}, {\em Modeling of polyatomic gases in the kinetic {{Fokker-Planck}} method by extension of the master equation approach}, AIP Conference Proceedings, 2996 (2024), p.~060004, \url{https://doi.org/10.1063/5.0187658}, \url{https://doi.org/10.1063/5.0187658} (accessed 2025-01-14).

\bibitem{bhatnagarModelCollisionProcesses1954}
{\sc P.~L. Bhatnagar, E.~P. Gross, and M.~Krook}, {\em A {{Model}} for {{Collision Processes}} in {{Gases}}. {{I}}. {{Small Amplitude Processes}} in {{Charged}} and {{Neutral One-Component Systems}}}, Physical Review, 94 (1954), pp.~511--525, \url{https://doi.org/10.1103/PhysRev.94.511}, \url{https://link.aps.org/doi/10.1103/PhysRev.94.511} (accessed 2024-06-25).

\bibitem{birdOneDimensionalSteadyFlows1994}
{\sc G.~A. Bird}, {\em One-{{Dimensional Steady Flows}}}, in Molecular {{Gas Dynamics And The Direct Simulation Of Gas Flows}}, G.~A. Bird, ed., Oxford University Press, May 1994, p.~257–315, \url{https://doi.org/10.1093/oso/9780198561958.003.0012}.

\bibitem{cercignaniChapterBoltzmannEquation2002}
{\sc C.~Cercignani}, {\em Chapter 1 - {{The Boltzmann Equation}} and {{Fluid Dynamics}}}, vol.~1 of Handbook of {{Mathematical Fluid Dynamics}}, North-Holland, 2002, pp.~1--69, \url{https://doi.org/10.1016/S1874-5792(02)80003-9}, \url{https://www.sciencedirect.com/science/article/pii/S1874579202800039}.

\bibitem{chaintronPropagationChaosReview2022}
{\sc L.-P. Chaintron and A.~Diez}, {\em Propagation of chaos: A review of models, methods and applications. {{I}}. {{Models}} and methods}, Kinetic and Related Models, 15 (2022), p.~895, \url{https://doi.org/10.3934/krm.2022017}, \url{https://arxiv.org/abs/2203.00446}.

\bibitem{chapmanMathematicalTheoryNonuniform1990}
{\sc S.~Chapman and T.~G. Cowling}, {\em The {{Mathematical Theory}} of {{Non-uniform Gases}}: {{An Account}} of the {{Kinetic Theory}} of {{Viscosity}}, {{Thermal Conduction}} and {{Diffusion}} in {{Gases}}}, Cambridge University Press, 1990.

\bibitem{OptimalPredictionMoria}
{\sc A.~J. Chorin, O.~H. Hald, and R.~Kupferman}, {\em Optimal prediction and the mori–zwanzig representation of irreversible processes}, Proceedings of the National Academy of Sciences, 97 (2000), pp.~2968--2973, \url{https://doi.org/10.1073/pnas.97.7.2968}, \url{https://www.pnas.org/doi/abs/10.1073/pnas.97.7.2968}, \url{https://arxiv.org/abs/https://www.pnas.org/doi/pdf/10.1073/pnas.97.7.2968}.

\bibitem{ewartMassFlowRate2007}
{\sc T.~Ewart, P.~Perrier, I.~Graur, and J.~M{\'e}olans}, {\em Mass flow rate measurements in a microchannel, from hydrodynamic to near free molecular regimes}, Journal of Fluid Mechanics, 584 (2007), pp.~337--356, \url{https://doi.org/10.1017/S0022112007006374}.

\bibitem{fisherTheoryStatisticalEstimation1925}
{\sc R.~A. Fisher}, {\em Theory of {{Statistical Estimation}}}, Mathematical Proceedings of the Cambridge Philosophical Society, 22 (1925), pp.~700--725, \url{https://doi.org/10.1017/S0305004100009580}, \url{https://www.cambridge.org/core/journals/mathematical-proceedings-of-the-cambridge-philosophical-society/article/abs/theory-of-statistical-estimation/7A05FB68C83B36C0E91D42C76AB177D4} (accessed 2024-06-28).

\bibitem{gorji2015variance}
{\sc M.~H. Gorji, N.~Andric, and P.~Jenny}, {\em Variance reduction for fokker--planck based particle monte carlo schemes}, Journal of Computational Physics, 295 (2015), pp.~644--664, \url{https://doi.org/10.1016/j.jcp.2015.04.008}.

\bibitem{gorjiEfficientParticleFokker2014}
{\sc M.~H. Gorji and P.~Jenny}, {\em An efficient particle {{Fokker}}--{{Planck}} algorithm for rarefied gas flows}, Journal of Computational Physics, 262 (2014), pp.~325--343, \url{https://doi.org/10.1016/j.jcp.2013.12.046}, \url{https://www.sciencedirect.com/science/article/pii/S0021999113008541} (accessed 2024-12-10).

\bibitem{gorjiFokkerPlanckDSMC2015}
{\sc M.~H. Gorji and P.~Jenny}, {\em Fokker--{{Planck}}--{{DSMC}} algorithm for simulations of rarefied gas flows}, Journal of Computational Physics, 287 (2015), pp.~110--129, \url{https://doi.org/10.1016/j.jcp.2015.01.041}.

\bibitem{gorjiEntropicFokkerPlanckKinetic2021}
{\sc M.~H. Gorji and M.~Torrilhon}, {\em Entropic {{Fokker-Planck}} kinetic model}, Journal of Computational Physics, 430 (2021), p.~110034, \url{https://doi.org/10.1016/j.jcp.2020.110034}, \url{https://www.sciencedirect.com/science/article/pii/S0021999120308081}.

\bibitem{gorjiFokkerPlanckModel2011}
{\sc M.~H. Gorji, M.~Torrilhon, and P.~Jenny}, {\em Fokker--{{Planck}} model for computational studies of monatomic rarefied gas flows}, Journal of Fluid Mechanics, 680 (2011), pp.~574--601, \url{https://doi.org/10.1017/jfm.2011.188}, \url{https://www.cambridge.org/core/product/identifier/S0022112011001881/type/journal_article} (accessed 2024-02-01).

\bibitem{gottwald2025stable}
{\sc G.~A. Gottwald, F.~Li, Y.~Marzouk, and S.~Reich}, {\em Stable generative modelling using schr{\"o}dinger bridges}, Philosophical Transactions A, 383 (2025), p.~20240332.

\bibitem{guptaAutomatedBoltzmannCollision2012}
{\sc V.~K. Gupta and M.~Torrilhon}, {\em Automated {{Boltzmann}} collision integrals for moment equations}, AIP Conference Proceedings, 1501 (2012), pp.~67--74, \url{https://doi.org/10.1063/1.4769474}.

\bibitem{heppKineticFokkerPlanck2020}
{\sc C.~Hepp, M.~Grabe, and K.~Hannemann}, {\em A kinetic {{Fokker}}--{{Planck}} approach to model hard-sphere gas mixtures}, Physics of Fluids, 32 (2020), p.~027103, \url{https://doi.org/10.1063/1.5141909}, \url{https://doi.org/10.1063/1.5141909} (accessed 2025-01-14).

\bibitem{jennySolutionAlgorithmFluid2010}
{\sc P.~Jenny, M.~Torrilhon, and S.~Heinz}, {\em A solution algorithm for the fluid dynamic equations based on a stochastic model for molecular motion}, Journal of Computational Physics, 229 (2010), pp.~1077--1098, \url{https://doi.org/10.1016/j.jcp.2009.10.008}, \url{https://linkinghub.elsevier.com/retrieve/pii/S0021999109005531} (accessed 2024-02-01).

\bibitem{khasminskiiStochasticStabilityDifferential2012}
{\sc R.~Khasminskii}, {\em Stochastic {{Stability}} of {{Differential Equations}}}, vol.~66 of Stochastic {{Modelling}} and {{Applied Probability}}, Springer, Berlin, Heidelberg, 2012, \url{https://doi.org/10.1007/978-3-642-23280-0}.

\bibitem{kirkwoodStatisticalMechanicalTheory1946}
{\sc J.~G. Kirkwood}, {\em The {{Statistical Mechanical Theory}} of {{Transport Processes I}}. {{General Theory}}}, The Journal of Chemical Physics, 14 (1946), pp.~180--201, \url{https://doi.org/10.1063/1.1724117}, \url{https://doi.org/10.1063/1.1724117} (accessed 2025-04-30).

\bibitem{knudsenGesetzeMolekularstromungUnd2006}
{\sc M.~Knudsen}, {\em Die {{Gesetze}} der {{Molekularstr{\"o}mung}} und der inneren {{Reibungsstr{\"o}mung}} der {{Gase}} durch {{R{\"o}hrn}}}, Annalen der Physik, 333 (2006), pp.~75--130, \url{https://doi.org/10.1002/andp.19093330106}.

\bibitem{pmlr-v202-lai23d}
{\sc C.-H. Lai, Y.~Takida, N.~Murata, T.~Uesaka, Y.~Mitsufuji, and S.~Ermon}, {\em {FP}-diffusion: Improving score-based diffusion models by enforcing the underlying score fokker-planck equation}, in Proceedings of the 40th International Conference on Machine Learning, A.~Krause, E.~Brunskill, K.~Cho, B.~Engelhardt, S.~Sabato, and J.~Scarlett, eds., vol.~202 of Proceedings of Machine Learning Research, PMLR, 23--29 Jul 2023, pp.~18365--18398, \url{https://proceedings.mlr.press/v202/lai23d.html}.

\bibitem{lebowitzNonequilibriumDistributionFunctions1960}
{\sc J.~L. Lebowitz, H.~L. Frisch, and E.~Helfand}, {\em Nonequilibrium {{Distribution Functions}} in a {{Fluid}}}, The Physics of Fluids, 3 (1960), pp.~325--338, \url{https://doi.org/10.1063/1.1706037}.

\bibitem{mathiaudFokkerPlanckModel2016}
{\sc J.~Mathiaud and L.~Mieussens}, {\em A {{Fokker}}--{{Planck Model}} of the {{Boltzmann Equation}} with {{Correct Prandtl Number}}}, Journal of Statistical Physics, 162 (2016), pp.~397--414, \url{https://doi.org/10.1007/s10955-015-1404-9}, \url{https://doi.org/10.1007/s10955-015-1404-9} (accessed 2025-02-01).

\bibitem{mathiaudFokkerPlanckModel2017}
{\sc J.~Mathiaud and L.~Mieussens}, {\em A {{Fokker}}--{{Planck Model}} of the {{Boltzmann Equation}} with {{Correct Prandtl Number}} for {{Polyatomic Gases}}}, Journal of Statistical Physics, 168 (2017), pp.~1031--1055, \url{https://doi.org/10.1007/s10955-017-1837-4}, \url{https://doi.org/10.1007/s10955-017-1837-4} (accessed 2025-01-14).

\bibitem{mckeanCLASSMARKOVPROCESSES1966}
{\sc H.~P. McKean}, {\em A {{Class of Markov Processes Associated with Nonlinear Parabolic Equations}}}, Proceedings of the National Academy of Sciences of the United States of America, 56 (1966), pp.~1907--1911.

\bibitem{pawulaApproximationLinearBoltzmann1967}
{\sc R.~F. Pawula}, {\em Approximation of the {{Linear Boltzmann Equation}} by the {{Fokker-Planck Equation}}}, Physical Review, 162 (1967), pp.~186--188, \url{https://doi.org/10.1103/PhysRev.162.186}, \url{https://link.aps.org/doi/10.1103/PhysRev.162.186} (accessed 2025-04-30).

\bibitem{pengStatisticalFluidMechanics2023}
{\sc H.~Peng}, {\em Statistical fluid mechanics: {{Dynamics}} equations and linear response theory}, Physics of Fluids, 35 (2023), p.~071704, \url{https://doi.org/10.1063/5.0156582}.

\bibitem{riskenFokkerPlanckEquationMethods1996}
{\sc H.~Risken}, {\em The {{Fokker-Planck Equation}}: {{Methods}} of {{Solution}} and {{Applications}}}, vol.~18 of Springer {{Series}} in {{Synergetics}}, Springer, Berlin, Heidelberg, 1996, \url{https://doi.org/10.1007/978-3-642-61544-3}, \url{https://link.springer.com/10.1007/978-3-642-61544-3} (accessed 2025-04-23).

\bibitem{schaafFlowRarefiedGases2015}
{\sc S.~A. Schaaf and P.~L. Chambre}, {\em H. {{Flow}} of {{Rarefied Gases}}}, in Fundamentals of {{Gas Dynamics}}, Princeton University Press, Dec. 2015, ch.~Fundamentals of Gas Dynamics, pp.~687--740.

\bibitem{SCHOULER2020100638}
{\sc M.~Schouler, Y.~Pr{\'e}vereaud, and L.~Mieussens}, {\em Survey of flight and numerical data of hypersonic rarefied flows encountered in earth orbit and atmospheric reentry}, Progress in Aerospace Sciences, 118 (2020), p.~100638, \url{https://doi.org/10.1016/j.paerosci.2020.100638}.

\bibitem{SHEN2003512}
{\sc C.~Shen, J.~Fan, and C.~Xie}, {\em Statistical simulation of rarefied gas flows in micro-channels}, Journal of Computational Physics, 189 (2003), pp.~512--526, \url{https://doi.org/10.1016/S0021-9991(03)00231-6}.

\bibitem{shermanAdjustmentInverseMatrix1950}
{\sc J.~Sherman and W.~J. Morrison}, {\em Adjustment of an {{Inverse Matrix Corresponding}} to a {{Change}} in {{One Element}} of a {{Given Matrix}}}, The Annals of Mathematical Statistics, 21 (1950), pp.~124--127, \url{https://doi.org/10.1214/aoms/1177729893}, \url{https://projecteuclid.org/journals/annals-of-mathematical-statistics/volume-21/issue-1/Adjustment-of-an-Inverse-Matrix-Corresponding-to-a-Change-in/10.1214/aoms/1177729893.full} (accessed 2025-03-31).

\bibitem{toscaniHtheoremAsymptoticTrend1987}
{\sc G.~Toscani}, {\em H-theorem and asymptotic trend of the solution for a rarefied gas in the vacuum}, Archive for Rational Mechanics and Analysis, 100 (1987), pp.~1--12, \url{https://doi.org/10.1007/BF00281245}, \url{https://doi.org/10.1007/BF00281245} (accessed 2025-07-04).

\bibitem{toscaniEntropyProductionRate1999}
{\sc G.~Toscani}, {\em Entropy {{Production}} and the {{Rate}} of {{Convergence}} to {{Equilibrium}} for the {{Fokker-Planck Equation}}}, Quarterly of Applied Mathematics, 57 (1999), pp.~521--541, \url{https://www.jstor.org/stable/43638304} (accessed 2025-01-07), \url{https://arxiv.org/abs/43638304}.

\bibitem{toscaniSharpEntropyDissipation1999}
{\sc G.~Toscani and C.~Villani}, {\em Sharp {{Entropy Dissipation Bounds}} and {{Explicit Rate}} of {{Trend}} to {{Equilibrium}} for the {{Spatially Homogeneous Boltzmann Equation}}}, Communications in Mathematical Physics, 203 (1999), pp.~667--706, \url{https://doi.org/10.1007/s002200050631}, \url{https://doi.org/10.1007/s002200050631} (accessed 2025-06-10).

\bibitem{truesdellTrendEquilibriumAccording1984a}
{\sc C.~Truesdell}, {\em The {{Trend}} to {{Equilibrium According}} to the {{Kinetic Theory}} of {{Gases}}}, Springer New York, New York, NY, 1984, pp.~451--459, \url{https://doi.org/10.1007/978-1-4612-5206-1_27}.

\bibitem{wangNumericalInvestigationFlow2024}
{\sc Q.~Wang, K.~Wang, X.~Wu, and Z.~Gao}, {\em Numerical investigation of flow and particles contamination in reticle mini environment for extreme ultraviolet lithography}, Journal of Vacuum Science \& Technology B, 42 (2024), p.~052602, \url{https://doi.org/10.1116/6.0003791}.

\bibitem{yunEntropyProductionEllipsoidal2016}
{\sc S.-B. Yun}, {\em Entropy production for ellipsoidal {{BGK}} model of the {{Boltzmann}} equation}, Kinetic and Related Models, 9 (Sun May 01 00:00:00 UTC 2016), pp.~605--619, \url{https://doi.org/10.3934/krm.2016009}, \url{https://www.aimsciences.org/en/article/doi/10.3934/krm.2016009} (accessed 2025-04-23).

\end{thebibliography}

\end{document}